\PassOptionsToPackage{dvipsnames}{xcolor}
\documentclass[letterpaper]{amsart}
\usepackage[margin = 2cm]{geometry}
\usepackage{amsaddr}

\usepackage{amsfonts}
\usepackage{amssymb}
\usepackage{tikz}
\usetikzlibrary{decorations.markings}
\usetikzlibrary{calc}
\usepackage[shortlabels]{enumitem}
\usepackage[dvipsnames]{xcolor}
\usepackage{graphicx}
\usepackage{subcaption}
\usepackage{float}
\labelformat{subfigure}{\thefigure(\textsc{#1})}
\usepackage[hidelinks]{hyperref}
\usepackage{url}
\usepackage{comment}

\usepackage[numbers]{natbib}

\usepackage{bm}

\usepackage{concmath}
\usepackage[slantedGreek]{ccfonts}

\usepackage{eulervm}

\usepackage{soul}

\theoremstyle{plain}
\newtheorem{thm}{Theorem}[section]

\newtheorem{lem}[thm]{Lemma}

\newtheorem{cor}[thm]{Corollary}
\newtheorem{conj}[thm]{Conjecture}

\newtheorem{obs}[thm]{Observation}

\theoremstyle{definition}
\newtheorem{defn}[thm]{Definition}

\theoremstyle{remark}
\newtheorem*{rmk}{Remark}

\def\deg{{\sf deg}}

\def\sqr#1#2{{\vcenter{\hrule height.#2pt
        \hbox{\vrule width.#2pt height#1pt \kern#1pt
                \vrule width.#2pt}
        \hrule height.#2pt}}}

\colorlet{niceblue}{blue!65!black!75}

\definecolor{mnt}{RGB}{167,230,215}
\definecolor{brwn}{RGB}{140, 70, 20}

\newcommand{\av}{\overline{\sigma}}

\hypersetup{
    colorlinks=true,  
    urlcolor=blue,
    }

\usepackage{todonotes}

\title{Proximity and Remoteness in Trees with given Segment Sequence}

\author{
Imran Allie\textsuperscript{1}, 
Brandon Du Preez\textsuperscript{1}, 
Zion Hefty\textsuperscript{4}, 
Sonwabile~Mafunda\textsuperscript{2,3},  
}

\address{
\textsuperscript{1}University~of~Cape~Town, 
\textsuperscript{2}University~of~Johannesburg,~South Africa,\\ 
\textsuperscript{3}Soka~University~of~America,
\textsuperscript{4}University~of~Denver,~USA
}

\email{imran.allie@uct.ac.za, brandon.dupreez@uct.ac.za, zion.hefty@du.edu, smafunda@soka.edu}

\begin{document}

\begin{abstract}
The average distance of a vertex $v$ of a connected graph $G$ is the arithmetic mean of the distances 
from $v$ to all other vertices of $G$. The proximity $\pi(G)$ and the remoteness $\rho(G)$ of $G$ 
are the minimum and the maximum of the average distances of the vertices of $G$, respectively. 
A segment $S$ of a tree is a path-subtree whose terminal vertices are either branching or pendant vertices, and all internal vertices have degree $2$. 
The segment sequence $(\ell_{S_1}, \ell_{S_2},\ldots,\ell_{S_m})$ of T is the non-increasing sequence of the lengths of all segments of T.

In this paper, we characterise trees with a given segment sequence that minimise proximity and those that maximise remoteness. For trees that minimise remoteness, we establish sufficient conditions under which the tree is starlike and conjecture that a remoteness-minimising tree is always either starlike or a ($\Delta,2$)-tree. For trees that maximise proximity, we derive a collection of necessary conditions that any tree maximising proximity must satisfy.
\end{abstract}

\maketitle
Keywords: Distance; proximity; remoteness; trees; segment sequence

MSC-class: 05C05; 05C12 \\

\footnotetext{This research supported in part by the NGA(MaSS) grant UCDP-712, CoE(MaSS) grant 2022-044-GRA-Workshop, NITheCS - Mathematical Structures research grant, an the University of Johannesburg GES-4IR grant.}


\section{Introduction}

Let $G$ be a connected graph of order $n$ with vertex set $V(G)$. The
\textit{average distance} $\overline{\sigma}(v,G)$ of a vertex $v$ of $G$
is defined as the arithmetic mean of the distances from $v$ to all other
vertices of $G$, that is,
\[
\overline{\sigma}(v,G)
=\frac{1}{n-1}\sum_{u\in V(G)}d(v,u),
\]
where $d(v,u)$ denotes the usual shortest-path distance between vertices
$v$ and $u$. The \textit{proximity} and the \textit{remoteness} of $G$,
denoted by $\pi(G)$ and $\rho(G)$, are the smallest and largest,
respectively, average distance among the vertices of $G$.

The average distance of a vertex is closely related to its
\textit{transmission} (also called its total distance), defined by $\sigma(v,G)=\sum_{u\in V(G)}d(v,u)$.
Clearly, $\sigma(v,G)=(n-1)\overline{\sigma}(v,G)$. Thus, extremal questions concerning proximity and remoteness may equivalently be formulated in terms of minimum and maximum transmission.

Bounds on proximity and remoteness in terms of order alone were given by
Zelinka~\cite{Zel1968} and later, independently, by Aouchiche and
Hansen~\cite{AouHan2011}, who introduced the names
\textit{proximity} and \textit{remoteness}. The following theorem describes
the extremal values when only the order of the graph is prescribed.

\begin{thm}[Zelinka~\cite{Zel1968}; Aouchiche and Hansen~\cite{AouHan2011}]
Let $G$ be a connected graph of order $n\geq 2$. Then
\[
1\leq \pi(G)\leq
\begin{cases}
\dfrac{n+1}{4}, & \text{if $n$ is odd},\\[6pt]
\dfrac{n+1}{4}+\dfrac{1}{4(n-1)}, & \text{if $n$ is even},
\end{cases}
\]
where the lower bound is attained if and only if $G$ has a vertex of
degree $n-1$, and the upper bound is attained if and only if $G$ is a path
or a cycle. Moreover,
\[
1\leq \rho(G)\leq \frac{n}{2},
\]
where the lower bound is attained if and only if $G$ is a complete graph,
and the upper bound is attained if and only if $G$ is a path.
\label{thm:order_bounds}
\end{thm}

Theorem~\ref{thm:order_bounds} provides the natural starting point for
the study of proximity and remoteness under additional restrictions.
Several results have considered restrictions involving minimum and
maximum degree, diameter, the number of pendant vertices, the number
of vertices of odd degree, the number of vertices of degree two, and
prescribed degree sequences; see, for example,
\cite{Dan2015, DMM2022, PZ2021, RB2014, GZ2023}.
Proximity and remoteness have also been studied for various restricted
graph classes, including trees, triangle-free graphs, $C_4$-free graphs,
graphs of prescribed girth, and outerplanar graphs, and related
graph classes; see
\cite{AH2017, BES1997, CLZ2021, CDOS2021, CDOS2022, DJM2021, DMM2026}.
Further references can be found in the survey~\cite{AR2024}.

In this paper we impose a different structural restriction by prescribing
the \textit{segment sequence} of a tree. A segment of a tree is a
path-subtree whose end-vertices are either pendant or branch vertices and
whose internal vertices have degree two. The \textit{segment sequence} of
a tree is the non-increasing sequence of the lengths of its segments.
The use of segments in the study of trees has a long history in connection
with the Wiener index; see the survey of Dobrynin, Entringer and
Gutman~\cite{DEG2001}. In particular, they recorded the basic
characterization that a sequence $(t_1,t_2,\ldots,t_m)$ of positive
integers is the segment sequence of an $n$-vertex tree if and only if
$t_1+t_2+\cdots+t_m+1=n$. 

Extremal problems for trees with a prescribed segment sequence have
subsequently been studied for several distance-based invariants. Lin and
Song~\cite{LS2016} characterized the trees that minimise the Wiener
index for a given segment sequence, while Andriantiana, Wagner and Wang
\cite{AWW2016} characterized the corresponding maximum
and showed that a maximising tree is a quasi-caterpillar. The same
setting has also been considered for the Steiner Wiener index by Zhang,
Wang and Zhang~\cite{ZWZ2019}, where starlike trees and
quasi-caterpillars arise as extremal structures. More recently, Allie and
Wade~\cite{AW2025} characterized the trees that minimise the
eccentricity index for a given segment sequence, with the unique
starlike tree attaining the minimum.

These results suggest a natural question for proximity and remoteness.
Although the Wiener index and its variants aggregate distance information
over sets of vertices, proximity and remoteness are determined by the
vertices with minimum and maximum average distance, respectively. Thus,
the arrangement of a prescribed collection of segments may affect these
parameters in ways that are not captured by the corresponding aggregate
distance measures.

The aim of this paper is to determine the extremal trees for proximity and
remoteness when the segment sequence is prescribed. More precisely, for a
given segment sequence, we characterise the trees that minimise proximity, 
and trees that maximise remoteness. 
We establish sufficient conditions and propose a conjecture for trees 
that minimise remoteness, and derive necessary conditions for trees 
that maximise proximity.
The results show that the arrangement of the prescribed segments around the
branch vertices is central to the extremal behaviour of these parameters.
In particular, some of the structures familiar from related extremal
problems reappear, while the structures that minimise or maximise
proximity and remoteness need not coincide.

The remainder of the paper is organised as follows. Section~2 introduces
the terminology, notation, and preliminary observations used throughout
the paper. Section~3 is devoted to trees with minimum proximity, while
Section~4 considers trees with minimum remoteness. Trees with maximum
remoteness are treated in Section~5, and trees with maximum proximity are
treated in Section~6. We conclude in Section~7 with a summary of our results
and some open questions.

\section{Preliminaries}

We use the following notation. 
Let $G=(V,E)$ be a finite connected graph.
The \textit{order} of $G$ is $|V|$, and will usually be denoted by $n$.
Denote by $n_k(G)$ the number of degree $k$ vertices in $G$.
A subset $X$ of $V$ is \textit{dominating} if every vertex of $V-X$ has a neighbour in $X$.

The \textit{proximity} $\pi(G)$ is the minimum average distance of a vertex of $G$, and the \textit{remoteness} $\rho(G)$ is the maximum average distance. 
A vertex $v$ is a \textit{median vertex} if it has minimum average distance ($\av(v) = \pi(G)$), and is an \textit{extreme vertex}, also called an \textit{outpost}, if it has maximum average distance ($\av(v) = \rho(G)$).

A \textit{subdivision} of $G$ is a graph formed from $G$ by replacing some of its edges by paths.
Two graphs $G$ and $H$ are \textit{homeomorphic} if they have a common subdivision.

A \textit{pendant vertex}, or \textit{leaf}, is a vertex of degree $1$, and a \textit{branch vertex} is a vertex of degree at least $3$. A \textit{bud vertex}, or simply a \textit{bud}, is a branch vertex incident with a pendant segment.

A \textit{segment} of a tree is a maximal path whose internal vertices have degree two and whose endpoints are pendant vertices or branch vertices.
The \textit{length} $\ell_i$ of a segment is the number of edges in the segment. 
A segment is \textit{pendant}, or a \textit{leaf segment}, if one of its end-vertices is a leaf; otherwise, it is a \textit{non-pendant segment}. 
The vertices at the ends of a segment are called its \textit{end-vertices}. 
Thus, an end-vertex is either a leaf or a branch vertex.

Suppose now that $T=(V,E)$ is a tree, and that $D = (\ell_1, \ell_2, \dots, \ell_k)$ is a non-increasing sequence of positive integers. 
We say that $T$ is a $D$\textit{-tree} if its segment sequence is $D$.
We use $S_i$ as the label of the segment which corresponds to length $l_i$.

The tree $T$ is a \textit{star} if it has a single dominating vertex, and is a \textit{caterpillar} if it has a dominating path. 
If $T$ has $h$ branch vertices and maximum degree $\Delta$, then we call it a \textit{$(\Delta,h)$-tree}.
A tree is \textit{starlike} if it is a $(\Delta,1)$-tree. 
If, in addition, $T$ is a $D$-tree, we call it a \textit{$(\Delta,h)$-$D$-tree}.

\begin{rmk}
    For a given segment sequence $D = (\ell_1, \ell_2, \dots, \ell_k)$, the $(k,1)$-$D$-tree is unique.
    \label{rem:unique_starlike}
\end{rmk}

A \textit{branch} of $T$ at a vertex $v$ is an induced subgraph of the form $T[C\cup \{v\}]$, where $C$ is a component of $T - \{v\}$.
A branch is \textit{trivial} if it consists of a single segment of $T$.
The \textit{branchweight} of a vertex $v$ of $T$, denoted by $bw(v)$, is the maximum size of a branch of $T$ at $v$.

The \textit{centroid} $C(T)$ of a tree $T$ is the set of vertices $v$ that minimize the branchweight $bw(v)$.
The following theorem of Jordan \cite{Jor1869} describes the structure of the centroid.
\begin{thm}[Jordan~\cite{Jor1869}]\label{thm:centroid_branchweight}
Let $T$ be a tree of order $n$. Then one of the following holds:
\begin{enumerate}[(i)]
\item $C(T)=\{c\}$ and $bw(c)\leq (n-1)/2$;
\item $C(T)=\{c_1,c_2\}$ and $bw(c_1)=bw(c_2)=n/2$.
\end{enumerate}
Moreover, if $v\in V(T)\setminus C(T)$, then $bw(v)> n/2$.
\end{thm}

Zelinka \cite{Zel1968} characterized the median vertices of a tree in terms of its centroid.
\begin{thm}[Zelinka~\cite{Zel1968}]\label{thm:zelinka_centroid}
The set of median vertices of a tree $T$ is its centroid $C(T)$.
\end{thm}

Combining these two results gives the following corollary, which is recorded by Barefoot, Entringer, and Székely~\cite{BES1997}.

\begin{cor}\label{cor:median_branchweight_bound}
Let $T$ be a tree of order $n$. If $v$ is a median vertex of $T$, then
\[
bw(v)\leq\frac{n}{2}.
\]
\end{cor}

The following classical theorem of Zelinka characterises median vertices in terms of branchweight.

\begin{thm}[Zelinka~\cite{Zel1968}]
Let $T$ be a tree. A vertex $v$ is a median vertex of $T$ if and only if
\[
bw(v)=\min_{u\in V(T)} bw(u).
\]
\label{thm:median_branchweight}
\end{thm}

We next record some elementary observations concerning distances in trees.
Suppose $u$ and $v$ are adjacent vertices of some tree, and let $\bm{T_u^{uv}}$ and $\bm{T_v^{uv}}$ denote the components of $T-uv$ containing $u$ and $v$, respectively.
Then for any vertex $w$ in $T_v^{uv}$, we have $d(u,w) = d(v,w) + 1$, and vice-versa. 

\begin{obs}\label{obs:edge_dist_comparison}
    Let $u,v$ be adjacent vertices in a tree $T$.
    Then $\sigma(u) = \sigma(v) - |T_u^{uv}| + |T_v^{uv}|$.
\end{obs}

\begin{obs}\label{obs:treeweight_flow}
    Let $T$ be a tree, and let $u,v$ be adjacent vertices of $T$.
    If $|T_u^{uv}| < |T_v^{uv}|$, then $\av(v) < \av(u)$. 
\end{obs}

We also note that every outpost of a tree is a leaf. Indeed, suppose that $u$ is an outpost that is not a leaf. 
Since $u$ has at least two neighbours, the components of $T-u$ have orders summing to $n-1$, and hence one of them, say $T_v$, has order at most $(n-1)/2$, where $v$ is the neighbour of $u$ in that component. 
Consequently,
\[
|T_v^{uv}|\leq \frac{n-1}{2}<n-\frac{n-1}{2}\leq |T_u^{uv}|.
\]
By Observation~\ref{obs:treeweight_flow}, with $u$ and $v$ interchanged, we obtain
\[
\av(u)<\av(v),
\]
contradicting the assumption that $u$ is an outpost. 
Thus we have:
\begin{obs}\label{obs:outpost_is_leaf}
Every outpost of a tree is a leaf.
\end{obs}

We \textit{remove} a non-pendant segment by removing all its internal vertices and edges while retaining its end-vertices. 
We \textit{remove} a pendant segment by removing all its vertices and edges except the end-vertex that is not a leaf. Suppose that $S_i$ is a non-pendant segment with end-vertices $u_i$ and $v_i$. 
After removing $S_i$, let $T_{u_i}^{S_i}$ and $T_{v_i}^{S_i}$ denote the components containing $u_i$ and $v_i$, respectively. 
The superscript distinguishes these components from those obtained by deleting an edge. We also define the following operation on a tree which will be useful.

\begin{defn}
    Let $S_i$ be a non-pendant segment of a tree $T$. Let $u$ and $v$ be the end vertices of $S_i$. The \textit{starring operation on $S_i$} is as follows: $S_i$ is removed, $u$ and $v$ are identified to $w$, a pendant segment which has the same length as $S_i$ is added adjacent to $w$ and preserves the label $S_i$.   
\end{defn}

It is clear that the starring operation preserves the segment sequence of the tree. If the context is clear, the starring operation may also be applied in reverse to restore the original tree.

\medskip

\section{Trees with minimum proximity}

In this section we determine the structure of the tree that minimises proximity among trees of a given segment sequence. We show that the unique tree attaining this minimum is the starlike tree in which the prescribed segments are all incident with a
single branch vertex. We first show that a tree with minimum proximity must
have maximum degree equal to the number of segments, and hence must be
starlike. Since the starlike tree with a given segment sequence is unique, this identifies the unique tree attaining the minimum proximity.
\vspace{2em}

\begin{lem}\label{lemprox1}
Let $D=(\ell_1,\ell_2,\ldots,\ell_k)$ be a non-increasing sequence of
positive integers, where $k\geq 3$. If $T$ is a $D$-tree with minimum
proximity, then
$\Delta(T)=k$.
\end{lem}
\begin{proof}
Let $T$ be a $D$-tree with minimum proximity. Let $B(T)$ denote the set of branch vertices of $T$, since
$k=1+\sum_{b\in B(T)}\left(\deg(b)-1\right)$,
we have $\Delta(T)\leq k$. Hence, it suffices to show that $\Delta(T)$ cannot be less than $k$. 
Suppose, to the contrary, that
$\Delta(T)\leq k-1$.
Then $T$ has at least two branch vertices, and hence it contains a
non-pendant segment. Let $S_j$ be such a segment, with end-vertices $u$
and $w$, and let $\ell_j$ be its length.
Let $v$ be a vertex of $T$ such that $\overline{\sigma}(v,T)=\pi(T)$.
We may assume without loss of generality that
$d_T(v,u)\leq d_T(v,w)$.
Let $X=T_{u}^{S_j}$ and $Y=T_{w}^{S_j}$ denote the two parts of $T$ determined by $S_j$ (see Figure \ref{figAB}). 

\begin{figure}[H]
\begin{center}
 \begin{tikzpicture}
      [scale=0.6,inner sep=1mm,
   vertex/.style={circle,thick,draw},
   thickedge/.style={line width=2pt},
       dedge/.style={
        double,
        double distance=2pt,
        postaction={decorate},
        decoration={
            markings,
        }
    }] 
    \begin{scope}
     \draw[yscale=1] (-2,2) circle (2.0)node {$X$};
      \node[vertex] (a1) at (0,2) [fill=black] {};
      \node [below] at (0.3,1.8) {$u$};
     \draw[dedge] (0.25,2) -- (7,2);
      \node[vertex] (a2) at (7,2) [fill=black] {};
            \node [below] at (6.8,1.8) {$w$};
          \draw[yscale=1] (9,2) circle (2.0) node {$Y$};
             \end{scope}
     \end{tikzpicture}
     \end{center}
     \caption{The tree $T$}
     \label{figAB}
\end{figure}
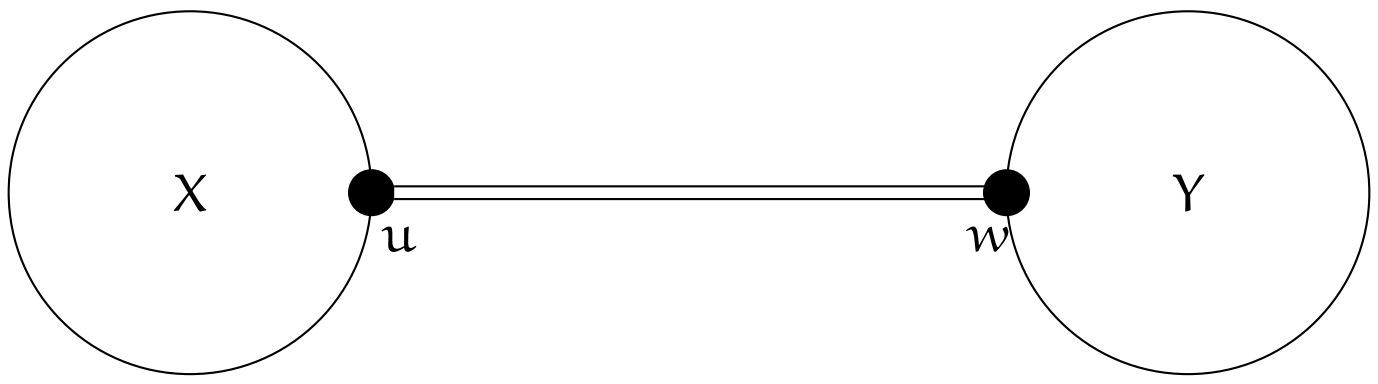

Since $w$ is a branch vertex and $S_j$ is
non-pendant, we have that $Y$ is non-trivial.
Construct a tree $T'$ from $T$ by moving all the edges incident with $w$
that do not belong to $S_j$ from $w$ to $u$. Equivalently, for every
vertex $z\in N_T(w)\setminus V(S_j)$, delete the edge $wz$ and add the
edge $uz$ (see Figure \ref{figAB'}). This is precisely the starring operation described above.

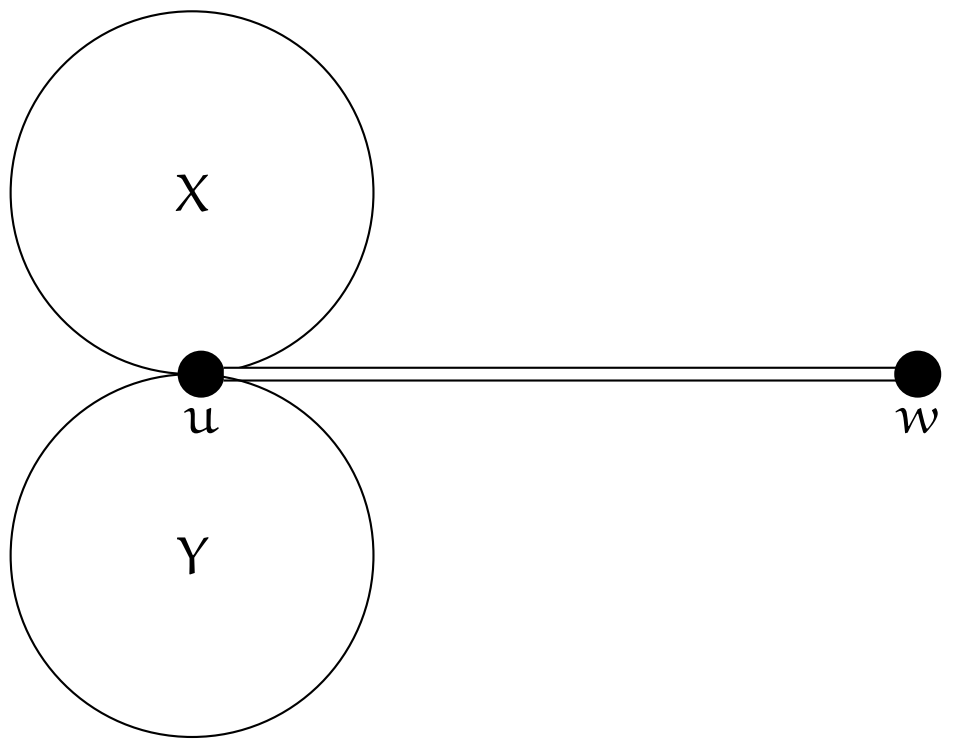
\begin{figure}[H]
\begin{center}
 \begin{tikzpicture}
      [scale=0.6,inner sep=1mm, 
   vertex/.style={circle,thick,draw},
   thickedge/.style={line width=2pt},
          dedge/.style={
        double,
        double distance=2pt,
        postaction={decorate},
        decoration={
            markings,
        }
    }]
    \begin{scope}
     \draw[yscale=1] (-3,2) circle (2) node {$X$};
      \node[vertex] (a1) at (-2.9,0.0) [fill=black] {};
      \node [below] at (-2.9,-0.2) {$u$};
     \draw[dedge] (-2.65,0.0) -- (5,0.0);
      \node[vertex] (a2) at (5,0.0) [fill=black] {};
            \node [below] at (5,-0.2) {$w$};
          \draw[yscale=1] (-3,-2) circle (2.0) node {$Y$};
             \end{scope}
     \end{tikzpicture}
     \end{center}
     \caption{The tree $T'$}
     \label{figAB'}
\end{figure}
\vspace{3em}

Consequently, $T'$ has the same segment sequence as $T$, and hence $T'$
is a $D$-tree.
For every $x\in X$ and every $y\in Y\setminus\{w\}$, the operation decreases
the distance from $x$ to $y$ by $\ell_j$. Thus
\[
\sigma(x, T')=\sigma(x, T)-\left(|Y|-1\right)\ell_j.
\]
If $v\in X$, it follows that $\sigma(v, T')<\sigma(v, T)$,
and therefore
$\pi(T')<\pi(T)$,
a contradiction.

\noindent It remains to consider the case
$v\in V(S_j)\setminus\{u,w\}$.
For every $y\in Y\setminus\{w\}$,
$d_T(v,y)=d_T(v,w)+d_T(w,y)$,
whereas
$d_{T'}(v,y)=d_T(v,u)+d_T(w,y)$.
Hence
\[
\sigma_{T'}(v)
=
\sigma_T(v)-\left(|Y|-1\right)
\bigl(d_T(v,w)-d_T(v,u)\bigr).
\]
By the choice of $u$, $d_T(v,u)\leq d_T(v,w)$,
so $\sigma_{T'}(v)\leq\sigma_T(v)$.
If the inequality is strict, then again
$\pi(T')<\pi(T)$,
a contradiction. We may therefore assume that
$d_T(v,u)=d_T(v,w)$.
Thus $v$ is the midpoint of $S_j$, and in particular $\ell_j$ is even.
Let $s$ be the neighbour of $v$ on the $u$-side of $S_j$. Then $d_T(s,w)-d_T(s,u)=2$.
Consequently, the starring operation gives
\[
\sigma(s, T')=\sigma(s, T)-2\left(|Y|-1\right).
\]
Let
$r=d_T(u,v)=\frac{\ell_j}{2}$. $T_s^{sv}$, the component of $T-sv$ containing $s$ has order
$|X|+r-1$. Since $v$ is a median vertex of $T$, Corollary~\ref{cor:median_branchweight_bound}
gives $bw(v)\leq \frac{n}{2}$. Hence, $|X|+r-1\leq \frac{n}{2}$.
Since
\[
n=|X|+|Y|+\ell_j-1=|X|+|Y|+2r-1,
\]
we obtain $|X|\leq |Y|+1$.
On the other hand, by Observation~\ref{obs:edge_dist_comparison},
\[
\sigma(s, T)-\sigma(v, T)
=
n-2(|X|+r-1)
=
|Y|-|X|+1.
\]
Since $|X|\geq 2$, we have $\sigma(s, T)-\sigma(v, T)\leq |Y|-1<2\left(|Y|-1\right)$.
Therefore,
\[
\begin{aligned}
\sigma(s, T')
&=\sigma(s, T)-2\left(|Y|-1\right)\\
&<\sigma(v, T).
\end{aligned}
\]
Since $\sigma(v, T')=\sigma(v, T)$ in the present case, we have
$\sigma(s, T')<\sigma(v, T')=\sigma(v, T)$.
Thus
\[
\pi(T')
\leq \frac{\sigma(s, T')}{n-1}
<
\frac{\sigma(v, T)}{n-1}
=\pi(T),
\]
again a contradiction.
Hence $\Delta(T)=k$.
\end{proof}

 Thus, any $D$-tree with minimum proximity must be starlike. Since the starlike $D$-tree is unique, this yields the following theorem. 
 
\begin{thm}\label{thm:min_proximity}
Let $D=\left(\ell_1, \ell_2, \ldots, \ell_k\right)$ be a segment sequence, where $k\geq 3$. Among all $D$-trees, the unique starlike $D$-tree has minimum proximity.
\end{thm}

\begin{proof}
There are finitely many $D$-trees, and hence there exists a $D$-tree $T$ having minimum proximity amongst all $D$-trees. By Lemma~\ref{lemprox1}, we have $\Delta(T)=k$. Thus, $T$ is starlike.
Since the starlike $D$-tree is unique, $T$ must be the unique starlike $D$-tree. Therefore, the unique starlike $D$-tree has minimum proximity amongst all $D$-trees.
 \end{proof}

\medskip


\section{Trees with minimum remoteness}

We have shown that for a segment sequence $D=(\ell_1, \ell_2, \ldots, \ell_k)$, the unique $(k,1)$-$D$-tree minimises proximity. The case to minimise remoteness turns out differently. In some cases, the unique $(k,1)$-tree minimises remoteness, but we have also found examples where the remoteness minimising tree is a $(\Delta,2)$-tree. For a segment sequence $D=(\ell_1, \ell_2, \ldots, \ell_k)$, we prove a sufficient condition on $D$ for the remoteness minimising $D$-tree to be a $(k,1)$-$D$-tree, and consequently as well for the remoteness minimising $D$-tree to be a $(\Delta,2)$-$D$-tree. We conjecture that the remoteness minimising $D$-tree is always a $(k,1)$-$D$-tree or a $(\Delta,2)$-$D$-tree. First, a useful lemma.

\begin{lem} \label{lem}
\label{lem:rmin_star_operation}
    Let $D=(\ell_1,\ell_2,\ldots,\ell_k)$ be a segment sequence and let $T$ be a $D$-tree. Let  $S_i$ be a leaf segment and $S_j$ be a non-leaf segment adjacent to $S_i$ such that $\ell_i \geq \ell_j$. Let $T'$ be the $D$-tree formed by performing the starring operation on $S_j$. Then $\rho(T') < \rho(T)$. 
\begin{proof}
    Let $T$ and $T'$ be as described in the statement of the Lemma. Let $u_j$ and $v_j$ be the end-vertices of segment $S_j$ in $T$, and also the leaf and non-leaf end-vertices of the leaf segment $S_j$ in $T'$, respectively. Let $X$ and $Y$ be the subtrees remaining after removal of segment $S_j$ from $T$. Let $u_i$ be the leaf vertex on segment $S_i$ in $T$ and $T'$. Then $\sigma(u_i,T') \geq \sigma(u_j,T')$ since $\ell_i \geq \ell_j$. We note the following: 

    \begin{itemize}
        \item[(1)] It is impossible for $u_j$ to be a unique outpost of $T'$. Moreover, if $u_j$ is an outpost of $T'$ then $\ell_i=\ell_j$, in which case, $u_i$ is also an outpost of $T'$.
        \item[(2)] For any vertex $x \in V(X) \cup V(Y)$, we have $\sigma(x,T') < \sigma(x,T)$.
    \end{itemize}
     
    \noindent \textit{Case 1.1:} $u_i$ is an outpost in $T$ and $T'$. 
        Then by (2) we have 
        $$\rho(T')=\sigma(u_i,T') < \sigma(u_i,T)=\rho(T).$$  
    \textit{Case 1.2:} $u_i$ is an outpost in $T$ and not in $T'$.
        Let $w'$ be an outpost in $T'$. By (1), $w' \neq u_j$ and obviously also $w' \notin S_j$. Therefore, by (2) we have 
        $$\rho(T')= \sigma(w,T') < \sigma(w,T) \leq \sigma(u_i,T) = \rho(T).$$  
    \textit{Case 2:} $u_i$ is not an outpost in $T$. 
        Let $w$ be an outpost in $T$. $w$ is a leaf and therefore $w \notin S_j$. Let $w'$ be an outpost in $T'$. Then 
        $$\rho(T')= \sigma(w',T') < \sigma(w',T) \leq \sigma(w,T) =\rho(T).$$   
\end{proof}
\end{lem}

If $D=(1,1,1,1,1)$, it is easy to determine that the remoteness minimising $D$-tree is the unique starlike tree. Let us instead consider the segment sequence $D=(10,1,1,1,1)$. We attempt to find the $D$-tree with the minimum remoteness. By Lemma \ref{lem}, no segment of length 1 can be a non-leaf segment, since it will necessarily be adjacent to some other segment which is at least the same size. Therefore, the only $D$-trees we need to consider is the $(5,1)$-$D$-tree, $T_1$, and the $(3,2)$-$D$-tree with the non-leaf segment of length 10, $T_2$. It turns out that $\rho(T_2) = 79 < \rho(T_1) = 85$. 

Given the above example, intuitively it seems that if segment lengths are sufficiently close ratio-wise, then the remoteness minimising tree is the unique starlike tree. While the existence of one segment sufficiently longer (ratio-wise) than any other segment would result in the starlike tree not being the one which minimises remoteness. This is due to the outpost changing from being a leaf vertex on the longest segment (being a large distance to all other leaf segments) to a leaf on one of the shorter segments (now not being a large distance to all leaf segments). In the following theorem we determine sufficient conditions on how much longer the long segment needs to be in order for the remoteness minimising tree to have two branching vertices.  

\begin{thm}
    Let $D=(\ell_1, \ell_2, \ldots, \ell_k)$ be a segment sequence with $\ell_1 = b\ell$ and $\ell_2=\ell_3=\dots=\ell_k=\ell$ and let $a=k-1$. If $a$ is even and $\frac{ab}{2} > a+b-1$, or $a$ is odd and $\frac{ab}{2} > \frac{2a+3b}{2}-1$, then the $D$-tree with minimum remoteness is a $(k,2)$-$D$-tree.
\begin{proof}
By Lemma \ref{lem}, it is clear that the remoteness minimising $D$-tree is either starlike or a $(\Delta,2)$-$D$-tree, since all segments of length $\ell$ would necessarily have to be leaf segments.

Let $T$ be the unique $(k,1)$-$D$-tree. Let $v$ be the leaf vertex on the segment of length $b\ell$. It is clear that $v$ is the unique outpost in $T$. Of all $(\Delta,2)$-$D$-trees, let $T'$ be the one which minimises remoteness. Obviously, segment $b\ell$ would be the non-leaf segment. Furthermore, the number of leaf segments on one end of the segment of length $b\ell$ would be $\lfloor \frac{a}{2} \rfloor$ and the number of leaf segments on the other end would be $\lceil \frac{a}{2} \rceil$. The outposts would be all the leaves on the former. Let one such outpost be $v$. See Figures \ref{fig_thm_1} and \ref{fig_thm_2}.  

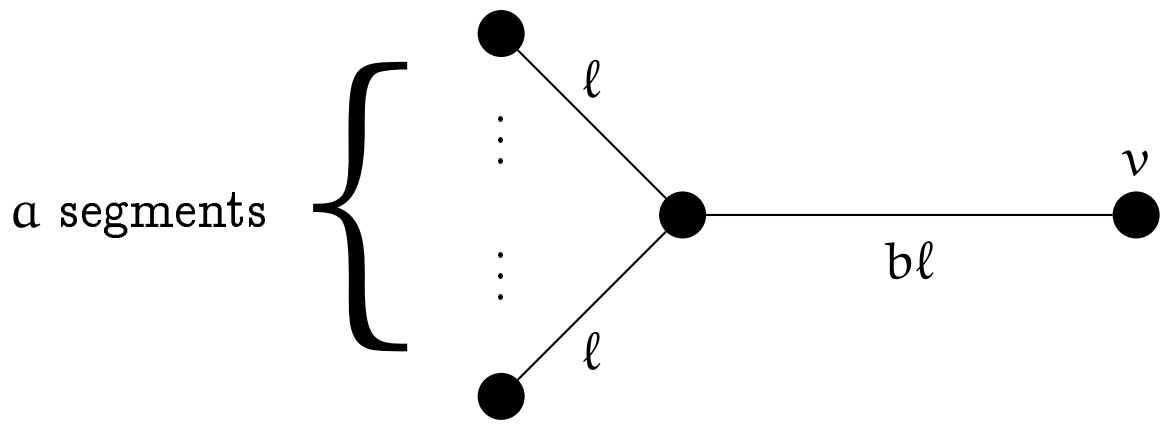
\begin{figure}[H]
\begin{center}
 \begin{tikzpicture}
      [scale=0.6,inner sep=1mm, 
   vertex/.style={circle,thick,draw}, 
   thickedge/.style={line width=2pt}] 
    \begin{scope}
        \node[vertex] (a1) at (0,0) [fill=black] {};
        \node[vertex] (a2) at (5,0) [fill=black] [label={$v$}] {};
        \node[vertex] (a3) at (-2,2) [fill=black] {};
        \node[vertex] (a4) at (-2,-2) [fill=black] {};
        
        \draw (a1) -- (a2);
        \node at (2.5,-0.5) {$b\ell$};
        \draw (a1) -- (a3);
        \node at (-1,1.5) {$\ell$};
        \draw (a1) -- (a4);
        \node at (-1,-1.5) {$\ell$};
        \node at (-2,1) {$\vdots$};
        \node at (-2,-0.5) {$\vdots$};

        \node at (-3.5,0) [scale=6] {$\{$};
        \node at (-6,0) {$a$ segments};
        
    \end{scope}
    \end{tikzpicture}
    \end{center}
    \caption{The unique starlike $D$-tree, $T$.}\label{fig_thm_1}
\end{figure}

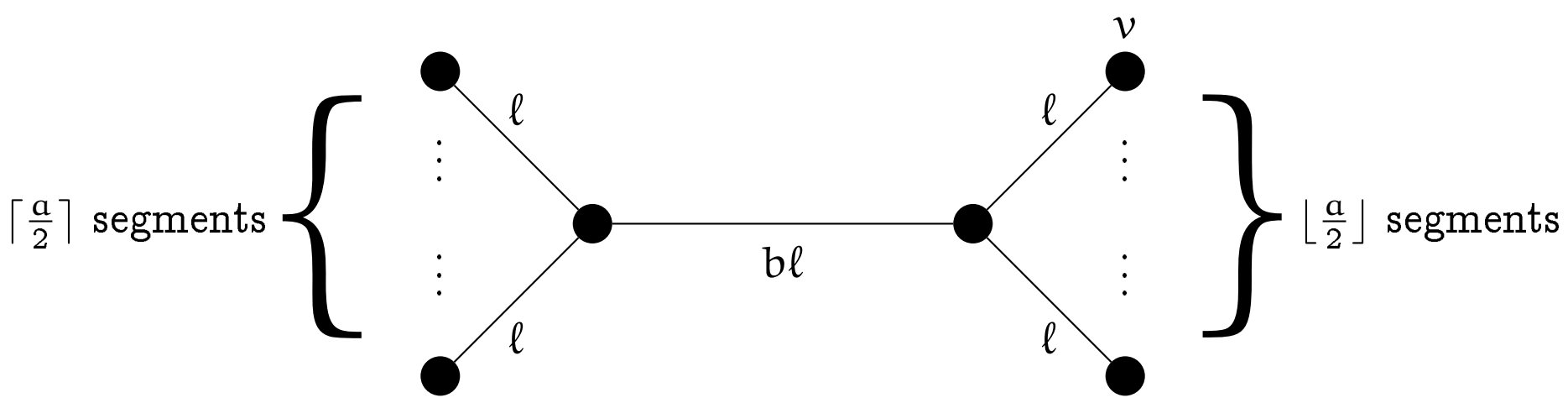
\begin{figure}[H]
\begin{center}
 \begin{tikzpicture}
      [scale=0.6,inner sep=1mm, 
   vertex/.style={circle,thick,draw}, 
   thickedge/.style={line width=2pt}] 
    \begin{scope}
        \node[vertex] (a1) at (0,0) [fill=black] {};
        \node[vertex] (a2) at (5,0) [fill=black] {};
        \node[vertex] (a3) at (-2,2) [fill=black] {};
        \node[vertex] (a4) at (-2,-2) [fill=black] {};
        \node[vertex] (a5) at (7,2) [fill=black,label={$v$}] {};
        \node[vertex] (a6) at (7,-2) [fill=black] {};
        
        \draw (a1) -- (a2);
        \node at (2.5,-0.5) {$b\ell$};
        \draw (a1) -- (a3);
        \node at (-1,1.5) {$\ell$};
        \draw (a1) -- (a4);
        \node at (-1,-1.5) {$\ell$};
        \draw (a2) -- (a5);
        \node at (6,1.5) {$\ell$};
        \draw (a2) -- (a6);
        \node at (6,-1.5) {$\ell$};

        \node at (-2,1) {$\vdots$};
        \node at (-2,-0.5) {$\vdots$};    
        \node at (7,1) {$\vdots$};
        \node at (7,-0.5) {$\vdots$};  

        \node at (-3.5,0) [scale=6] {$\{$};
        \node at (-6,0) {$\lceil \frac{a}{2} \rceil$ segments};

        \node at (8.5,0) [scale=6] {$\}$};
        \node at (11,0) {$\lfloor \frac{a}{2} \rfloor$ segments};
        
    \end{scope}
    \end{tikzpicture}
    \end{center}
    \caption{The remoteness minimising $(\Delta,2)$-$D$-tree, $T'$.} \label{fig_thm_2}
\end{figure}

Summing the distances from $v$ to each other vertex then simplifying, we have the following:
\begin{align*}
    \sigma(v,T) & = \frac{b\ell(b\ell+1)}{2} + \frac{a\ell(\ell+1)}{2} + ab\ell^2 \\
    \sigma(v,T') & = \lfloor \frac{a}{2} \rfloor \frac{\ell(\ell+1)}{2} +
                    (\lfloor \frac{a}{2} \rfloor - 1)\ell^2 + 
                    \frac{b\ell(b\ell+1)}{2} + b\ell^2 \\
                    & +
                    \lceil \frac{a}{2} \rceil \frac{\ell(\ell+1)}{2} + 
                    \lceil \frac{a}{2} \rceil (b+1)\ell^2
\end{align*}

If $a$ is even, then the inequality $\sigma(v,T) > \sigma(v,T')$ simplifies to $\frac{ab}{2} > a+b-1$. If $a$ is odd, then the inequality $\sigma(v,T) > \sigma(v,T')$ simplifies to $\frac{ab}{2} > \frac{2a+3b}{2}-1$.

\end{proof}
\end{thm}

\begin{cor}
    Let $D=(\ell_1, \ell_2, \ldots, \ell_k)$ be a segment sequence with $\ell_1 = b\ell$ and $\ell_2=\ell_3=\dots=\ell_k=\ell$ and let $a=k-1$. If $a$ is even and $\frac{ab}{2} \leq a+b-1$, or $a$ is odd and $\frac{ab}{2} \leq \frac{2a+3b}{2}-1$, then the $D$-tree with minimum remoteness is starlike.
\end{cor}

Given that having one segment sufficiently longer than the rest ratio-wise results in a remoteness minimising tree having two branching vertices, it is natural to consider the case where the longest segment is significantly longer than the second longest, which is then significantly longer than the rest. The idea being that should a segment sequence exist with a $(\Delta,3)$-tree being remoteness minimising, then it would likely have this form. With the aide of a computer we considered the following form of a segment sequence: $D=(\ell_1, \ell_2, \ldots, \ell_k)$ where $\ell_3=\ell_4=\dots=\ell_k=\ell$, $\ell_{2}=c\ell$ and $\ell_{1}=b\ell$, with $c < b$. By testing various values of $\ell,b,c$, with $b$ being much greater than $c$ and $c$ being much greater than $\ell$, we were not able to find such a segment sequence with a $(\Delta,3)$-tree being remoteness minimising. Therefore, we conjecture the following.

\medskip

\begin{conj}
    Let $D=(\ell_1,\ell_2,\ldots,\ell_k)$ be a segment sequence. Then the $D$-tree with minimum remoteness is starlike or a $(\Delta,2)$-$D$-tree. 
\end{conj}


\section{Trees with maximum remoteness}

In this section, we characterise the trees with fixed segment sequence that maximise remoteness. 
Let $D = (\ell_1, \ell_2, \dots, \ell_k)$ be a segment sequence. 
For $k$ odd, we construct the subdivided caterpillar $R_o(D)$ as follows. 
Begin with a caterpillar having $k$ edges, no degree 2 vertices and maximum degree 3.
Label the vertices of the dominating path $P$ by $P = u_0, u_1, \dots, u_{\frac{k+1}{2}}$. 
Then label the leaf adjacent to $u_1$ as $u_k$, the leaf adjacent to $u_2$ as $u_{k-1}$ and so on.
To obtain $R_o(D)$ from this caterpillar, replace the each edge $u_iu_j$, $i<j$, by a path of length $\ell_j$.

For $k$ even, we construct the caterpillar $R_e(D)$.
The construction is the same as that of $R_o(D)$, except that vertex $u_{\frac{k}{2} - 1}$ (second last vertex on the dominating path) has degree $4$: it is adjacent to $u_{\frac{k}{2} - 2}$, $u_{\frac{k}{2}}$, $u_{\frac{k}{2} + 1}$ and $u_{\frac{k}{2} + 2}$.

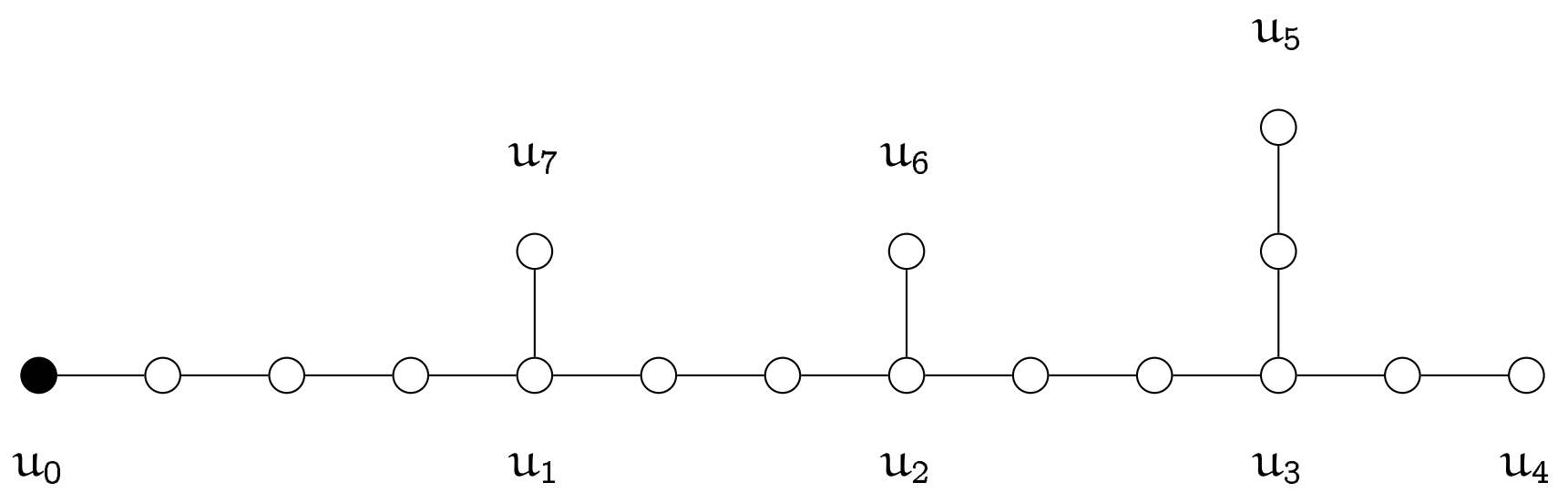
\begin{figure}[H]
    \centering
   \begin{tikzpicture}
    [scale=0.9,inner sep=0.9mm]
    \pgfdeclarelayer{nodelayer} 
    \pgfdeclarelayer{edgelayer}  
    \pgfsetlayers{edgelayer, nodelayer} 
    \tikzstyle{none}=[]
    \tikzstyle{whitevertex}=[fill=white, draw=black, shape=circle]
    \tikzstyle{blackvertex}=[fill=black, draw=black, shape=circle]
 
	\begin{pgfonlayer}{nodelayer}
		\node [blackvertex] (0) at (-6, 0) {};
		\node [whitevertex] (1) at (-5, 0) {};
		\node [whitevertex] (2) at (-4, 0) {};
		\node [whitevertex] (3) at (-3, 0) {};
		\node [whitevertex] (4) at (-2, 0) {};
		\node [whitevertex] (5) at (-1, 0) {};
		\node [whitevertex] (6) at (0, 0) {};
		\node [whitevertex] (7) at (1, 0) {};
		\node [whitevertex] (8) at (2, 0) {};
		\node [whitevertex] (9) at (3, 0) {};
		\node [whitevertex] (10) at (4, 0) {};
		\node [whitevertex] (11) at (5, 0) {};
		\node [whitevertex] (12) at (6, 0) {};
		\node [whitevertex] (13) at (-2, 1) {};
		\node [whitevertex] (14) at (1, 1) {};
		\node [whitevertex] (15) at (4, 1) {};
		\node [whitevertex] (16) at (4, 2) {};
		\node [] (17) at (-6, -0.75) {$u_0$};
		\node [] (18) at (-2, -0.75) {$u_1$};
		\node [] (19) at (1, -0.75) {$u_2$};
		\node [] (20) at (4, -0.75) {$u_3$};
		\node [] (21) at (6, -0.75) {$u_4$};
		\node [] (22) at (4, 2.75) {$u_5$};
		\node [] (23) at (1, 1.75) {$u_6$};
		\node [] (24) at (-2, 1.75) {$u_7$};
		\node [] (25) at (0, -1.75) {$R_o(4,3,3,2,2,1,1)$};
	\end{pgfonlayer}
	\begin{pgfonlayer}{edgelayer}
		\draw (0) to (1);
		\draw (1) to (2);
		\draw (2) to (3);
		\draw (3) to (4);
		\draw (4) to (5);
		\draw (5) to (6);
		\draw (6) to (7);
		\draw (7) to (8);
		\draw (8) to (9);
		\draw (9) to (10);
		\draw (10) to (11);
		\draw (11) to (12);
		\draw (13) to (4);
		\draw (14) to (7);
		\draw (16) to (15);
		\draw (15) to (10);
	\end{pgfonlayer}
\end{tikzpicture}
    \caption{The tree $R_o(4,3,3,2,2,1,1)$ maximises remoteness among all trees with segment sequence $(4,3,3,2,2,1,1)$.}
    \label{fig:remoteness_odd_extremal}
\end{figure}

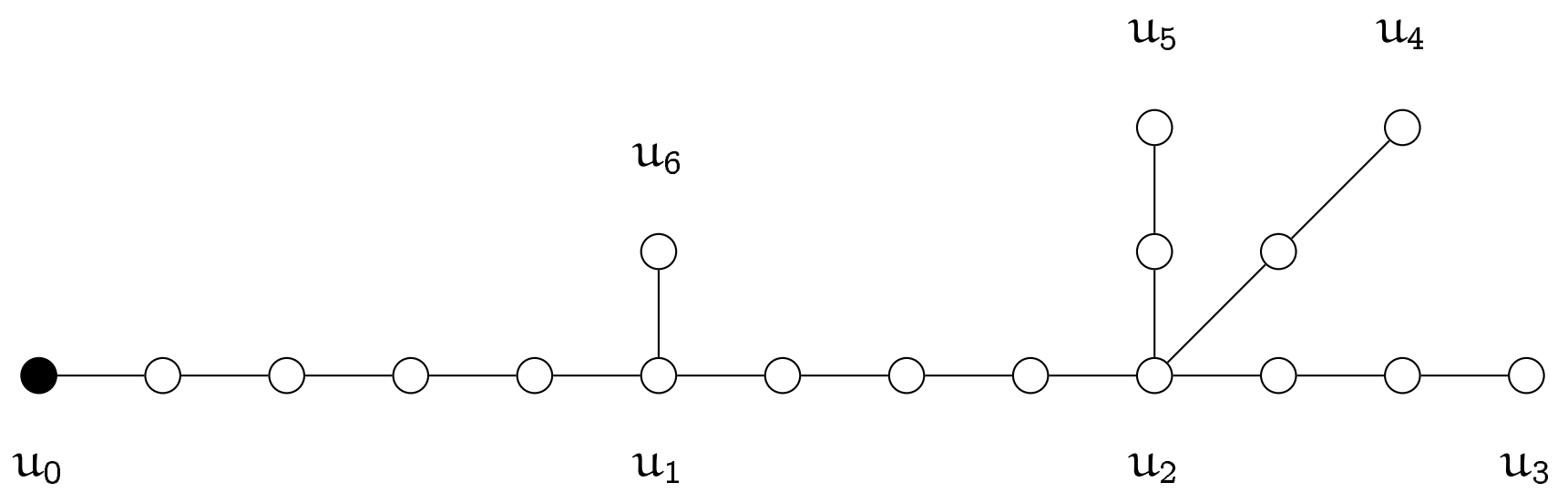
\begin{figure}[H]
    \centering
    \begin{tikzpicture}
    [scale=0.9,inner sep=0.9mm]
    \pgfdeclarelayer{nodelayer} 
    \pgfdeclarelayer{edgelayer}  
    \pgfsetlayers{edgelayer, nodelayer} 

    \tikzstyle{none}=[]
    \tikzstyle{whitevertex}=[fill=white, draw=black, shape=circle]
    \tikzstyle{blackvertex}=[fill=black, draw=black, shape=circle]
    
	\begin{pgfonlayer}{nodelayer}
		\node [whitevertex] (0) at (0, 0) {};
		\node [whitevertex] (1) at (-1, 0) {};
		\node [whitevertex] (2) at (1, 0) {};
		\node [whitevertex] (3) at (-2, 0) {};
		\node [whitevertex] (4) at (2, 0) {};
		\node [whitevertex] (5) at (-3, 0) {};
		\node [whitevertex] (6) at (3, 0) {};
		\node [whitevertex] (7) at (-4, 0) {};
		\node [whitevertex] (8) at (4, 0) {};
		\node [whitevertex] (9) at (-5, 0) {};
		\node [whitevertex] (10) at (5, 0) {};
		\node [blackvertex] (11) at (-6, 0) {};
		\node [whitevertex] (12) at (6, 0) {};
		\node [whitevertex] (13) at (4, 1) {};
		\node [whitevertex] (14) at (5, 2) {};
		\node [whitevertex] (16) at (3, 1) {};
		\node [whitevertex] (17) at (3, 2) {};
		\node [whitevertex] (18) at (-1, 1) {};
		\node [] (20) at (-1, 1.75) {$u_6$};
		\node [] (21) at (3, 2.75) {$u_5$};
		\node [] (22) at (5, 2.75) {$u_4$};
		\node [] (23) at (6, -0.75) {$u_3$};
		\node [] (24) at (3, -0.75) {$u_2$};
		\node [] (25) at (-1, -0.75) {$u_1$};
		\node [] (26) at (-6, -0.75) {$u_0$};
		\node [] (27) at (0, -1.75) {$R_e(5,4,3,2,2,1)$};
	\end{pgfonlayer}
	\begin{pgfonlayer}{edgelayer}
		\draw (11) to (9);
		\draw (9) to (7);
		\draw (7) to (5);
		\draw (5) to (3);
		\draw (3) to (1);
		\draw (1) to (0);
		\draw (0) to (2);
		\draw (2) to (4);
		\draw (4) to (6);
		\draw (6) to (8);
		\draw (8) to (10);
		\draw (10) to (12);
		\draw (6) to (13);
		\draw (13) to (14);
		\draw (6) to (16);
		\draw (16) to (17);
		\draw (1) to (18);
	\end{pgfonlayer}
\end{tikzpicture}
    \caption{The tree $R_e(5,3,3,2,2,1)$ maximises remoteness among all trees with segment sequence $(5,3,3,2,2,1)$.}
    \label{fig:remoteness_even_extremal}
\end{figure}
\medskip

\begin{lem}
    Let $D$ be a non-increasing sequence of $k$ positive integers.
    Then any $D$-tree $T$ maximising $\rho(T)$ is a subdivided caterpillar with the following properties:
    \begin{itemize}
        \item If $k$ is odd, then $\Delta(T) = 3$,
        \item If $k$ is even, then $\Delta(T) = 4$, and the unique vertex of degree 4 is incident with three pendant segments,
        \item For each outpost of $T$, there is a path containing both the outpost and every branch vertex.
    \end{itemize}
    \label{lem:remoteness_max}
\end{lem}

\begin{proof}

Let $T$ be a $D$-tree maximising $\rho(T)$.

\textit{Claim 1:} $\Delta(T) \leq 4$.\\
\textit{Proof of Claim 1:} Let $w$ be an outpost, $P: w, \dots, x$ a longest path starting at $w$, and assume to the contrary that there exists a vertex $v$ that has degree at least 5. 
Pull two subtrees off $v$, attach them to the end of the path $x$, thus creating a new tree $T'$ with the same segment sequence such that $\av(w, T') > \av(w, T)$ (and thus $\rho(T') > \rho(T)$).

\textit{Claim 2:} $T$ is a subdivided caterpillar. Further, for each outpost, there is a path in the caterpillar that contains the outpost and every branch vertex.\\
\textit{Proof of Claim 2: } Let $w$ be an outpost, and $P: w,\dots, x$ a longest path starting at $w$.
Assume to the contrary that $P$ does not contain each branch vertex.
Consider the forest $T-E(P)$ formed by removing every edge of $P$. 
This forest contains some component $Y$ with a branch vertex.
There is a segment $S$ of $Y$ incident with $P$.
Let $Y'$ be the subgraph of $Y$ induced by the edges $E(Y) - E(S)$ (i.e., $Y'$ consists of all segments of $Y$, except for $S$).
Moving the sub-tree $Y'$ onto $x$ yields a tree $T'$ such that $\av(w, T') > \av(w,T)$.

\textit{Claim 3:} $T$ has at most one vertex of degree $4$.\\
\textit{Proof of Claim 3:} Suppose to the contrary that $T$ has two vertices of degree 4, $u$ and $v$.
Again let $w$ be an outpost of $T$, and $P = w,\dots, x$ be a longest path starting at $w$.
Let $S_u$ be a segment that is incident with $u$, but not in $P$; and $S_v$ a segment incident with $v$ but not in $P$. 
Moving both $S_u$ and $S_v$ onto $x$ yields a new tree $T'$ such that $\av(w, T') > \av(w,T)$.

\textit{Claim 4:} If $T$ has a vertex of degree $4$, then that vertex is incident with $3$ pendant segments.\\

\textit{Proof of Claim 4:} Assume the claim is false, and let $w$, $P$, $x$ all be defined as they are in Claim 3.
Let $u$ be the vertex of degree 4, and let $v$ be the branch vertex of $P$ that is furthest from $w$ (so $v$ is incident with two pendant segments, one of which contains $x$). 
By the previous claim, $u$ must lie on $P$.
Create a new $D$-tree $T$ by moving a segment of $T-P$ away from $u$, and make it incident with $v$ instead. 
Since $d(w,u) < d(w,v)$, we have $\av(w, T') > \av(w,T)$, completing the proof.
\end{proof}

\begin{thm}
    For a given segment sequence $D$, let $R(D)$ be $R_o(D)$ if $D$ has an odd number of entries and $R_e(D)$ if $D$ has an even number of entries. Then $R(D)$ is the unique (up to isomorphism) $D$-tree maximising remoteness.
    \label{thm:remoteness_minimiser}
\end{thm}

\begin{proof}
    Let $D = (\ell_1, \ell_2, \dots, \ell_k)$, $k\geq 3$ be a segment sequence. 
    Assume to the contrary that $T$ is a $D$-tree with $\rho(T) \geq \rho(R(D))$.
    Among all such $D$-trees, choose $T$ so that $\rho(T)$ attains the maximum possible value.
    Let $w$ be an outpost of $T$.
    Per Lemma \ref{lem:remoteness_max}, we may assume that $T$ is a subdivided caterpillar.
    Further, there is a maximal path $P : w,\dots, x$ of $T$ containing every branch vertex.
    Let $b_1, b_2, \dots, b_j$ be the sequence of these branch vertices, ordered as they are in $P$.
    Let $S_1$ be the segment starting at $w$ and ending at $b_1$, and let $S_i$ be the segment between $b_{i-1}$ and $b_i$ when $1<i\leq j$.
    Further, let $S_k$ be the so-far-unlabeled pendant segment incident with $b_1$, $S_{k-1}$ the pendant segment incident with $b_2$, and so on. 

    \textit{Claim 1:} Segment $S_i$ has length $\ell_i$ for $1\leq i \leq j$.\\
    Assume otherwise for the sake of contradiction, and let $q$ be the smallest integer such that $\ell(S_q) \neq \ell_q$. 
    There is a segment $S'$ of length $\ell_q$ somewhere in T. 
    By the minimality of $q$, at least one such segment $S'$ does not appear on the part of path $P$ from $w$ to the end of $S_q$ (i.e., $S'$ is either a pendant segment other than $S_1$, or is a segment $S_t$ for some $t>q$). 
    Note that $\ell(S') = \ell_q > \ell(S_q)$. 
    Create a new tree $T'$ by swapping $S_q$ and $S'$. 
    It is clear that $\av(w, T') > \av(w,T)$, contradicting the maximality of $\rho(T)$.

    \textit{Claim 2:} The segment $S_i$ has length $\ell_i$ for all $j < i \leq k$ (up to relabelling of the pendant segments incident with $b_j$.\\
    Assume otherwise for the sake of contradiction, and let $q$ be the largest integer such that $\ell(S_q) \neq \ell_q$.
    Some segment $S_t$ has length $\ell_q$. 
    Per Claim 1, $j < t < q$ (so $S_t$ is a pendant segment).
    By maximality of $q$, we have that $\ell(S_q) > \ell_q$.
    Thus we create a new $D$-tree $T'$ by swapping the segments $S_q$ and $S_t$, which yields $\av(w, T') > \av(w,T)$.
\end{proof}

\medskip


\section{Trees with maximum proximity}

Let $D$ be a sequence of non-increasing positive integers. 
We say that a $D$-tree $T$ \textbf{maximizes proximity} (or is $\pi$-maximizing) if $\pi(T) \geq \pi(T')$ for all $D$-trees $T'$.
In this section, we give a detailed set of necessary conditions that a $\pi$-maximising $D$-tree satisfies. 

\begin{thm}
    \label{thm:pmax_main_theorem}
    Let $D$ be a segment sequence of length $|D| \geq 5$, let $T$ be a $\pi$-maximizing $D$-tree and suppose $w$ is a median vertex.
    Then $T$ is a subdivided caterpillar with $\Delta(T) \in \{3,4\}$, and $w$ lies in a longest segment (possibly as an end-vertex). 
    There are (at most) two non-trivial branches at $w$, and each branch $B$ satisfies all the following:
    \begin{itemize}
        \item $B$ either has no degree 4 vertices, or it has one degree 4 vertex. In the second case, this vertex is the branch vertex of $B$ that is furthest from $w$. 
        \item Any leaf segment of $T$ in $B$ has length less than or equal to any non-leaf segment of $T$ in $B$. 
        \item If $X, Y$ are leaf segments such that $|X| > |Y|$, then $d(w,X) \geq d(w, Y)$. 
        \item If $X, Y$ are non-leaf segments such that $|X| > |Y|$, then $d(w,X) < d(w, Y)$. 
    \end{itemize}
\end{thm}

\begin{figure}[ht]
    \centering
    \begin{tikzpicture}[scale=0.8, inner sep = 0.85mm]
        \tikzstyle{none}=[]
        \tikzstyle{whitevertex}=[fill=white, draw=black, shape=circle]
        \tikzstyle{blackvertex}=[fill=black, draw=black, shape=circle]
		\node [style=whitevertex] (0) at (-6, 0) {};
		\node [style=whitevertex] (1) at (-5, 0) {};
		\node [style=whitevertex] (2) at (-3, 0) {};
		\node [style=whitevertex] (3) at (-5, 1) {};
		\node [style=whitevertex] (4) at (-3, 1) {};
		\node [style=whitevertex] (5) at (0, 0) {};
		\node [style=whitevertex] (6) at (0, 1) {};
		\node [style=whitevertex] (7) at (2, 0) {};
		\node [style=whitevertex] (8) at (2, 1) {};
		\node [style=whitevertex] (9) at (3.5, 0) {};
		\node [style=whitevertex] (10) at (3.5, 1) {};
		\node [style=whitevertex] (11) at (4.5, 0) {};
		\node [style=none] (12) at (-4, 0.3) {11};
		\node [style=none] (13) at (-1.5, 0.3) {25};
		\node [style=none] (14) at (1, 0.3) {10};
		\node [style=none] (15) at (2.8, 0.3) {4};
		\node [style=whitevertex] (16) at (-6, -2.5) {};
		\node [style=whitevertex] (17) at (-5, -2.5) {};
		\node [style=whitevertex] (18) at (-3, -2.5) {};
		\node [style=whitevertex] (19) at (-5, -1.5) {};
		\node [style=whitevertex] (20) at (-3, -1.5) {};
		\node [style=whitevertex] (21) at (0, -2.5) {};
		\node [style=whitevertex] (22) at (0, -1.5) {};
		\node [style=whitevertex] (23) at (2, -2.5) {};
		\node [style=whitevertex] (24) at (2, -1.5) {};
		\node [style=whitevertex] (25) at (3.5, -2.5) {};
		\node [style=whitevertex] (26) at (3.5, -1.5) {};
		\node [style=whitevertex] (27) at (4.5, -2.5) {};
		\node [style=none] (28) at (-4, -2.2) {10};
		\node [style=none] (29) at (-1.5, -2.2) {25};
		\node [style=none] (30) at (1, -2.2) {11};
		\node [style=none] (31) at (2.8, -2.2) {4};
		\node [style=none] (32) at (-7, 0.5) {$T_1$};
		\node [style=none] (33) at (-7, -2) {$T_2$};
		\node [style=whitevertex] (34) at (-6, -5) {};
		\node [style=whitevertex] (35) at (-4.5, -5) {};
		\node [style=whitevertex] (36) at (-2.5, -5) {};
		\node [style=whitevertex] (37) at (-4.5, -4) {};
		\node [style=whitevertex] (38) at (-2.5, -4) {};
		\node [style=whitevertex] (39) at (0.5, -5) {};
		\node [style=whitevertex] (40) at (0.5, -4) {};
		\node [style=whitevertex] (41) at (2.5, -5) {};
		\node [style=whitevertex] (42) at (2.5, -4) {};
		\node [style=whitevertex] (43) at (3.5, -5) {};
		\node [style=whitevertex] (44) at (3.5, -4) {};
		\node [style=whitevertex] (45) at (4.5, -5) {};
		\node [style=none] (46) at (-5.2, -4.7) {4};
		\node [style=none] (47) at (-1, -4.7) {25};
		\node [style=none] (48) at (-3.5, -4.7) {11};
		\node [style=none] (49) at (1.55, -4.7) {10};
		\node [style=none] (50) at (-7, -4.5) {$T_3$};
		\node [style=whitevertex] (51) at (-6, -7.5) {};
		\node [style=whitevertex] (52) at (-4.5, -7.5) {};
		\node [style=whitevertex] (53) at (-2.5, -7.5) {};
		\node [style=whitevertex] (54) at (-4.5, -6.5) {};
		\node [style=whitevertex] (55) at (-2.5, -6.5) {};
		\node [style=whitevertex] (56) at (0.5, -7.5) {};
		\node [style=whitevertex] (57) at (0.5, -6.5) {};
		\node [style=whitevertex] (58) at (2.5, -7.5) {};
		\node [style=whitevertex] (59) at (2.5, -6.5) {};
		\node [style=whitevertex] (60) at (3.5, -7.5) {};
		\node [style=whitevertex] (61) at (3.5, -6.5) {};
		\node [style=whitevertex] (62) at (4.5, -7.5) {};
		\node [style=none] (63) at (-5.2, -7.2) {4};
		\node [style=none] (64) at (-1, -7.2) {25};
		\node [style=none] (65) at (-3.5, -7.2) {10};
		\node [style=none] (66) at (1.55, -7.2) {11};
		\node [style=none] (67) at (-7, -7) {$T_4$};

		\draw (0) to (1);
		\draw (3) to (1);
		\draw (1) to (2);
		\draw (4) to (2);
		\draw (2) to (5);
		\draw (5) to (6);
		\draw (5) to (7);
		\draw (7) to (8);
		\draw (7) to (9);
		\draw (9) to (10);
		\draw (9) to (11);
		\draw (16) to (17);
		\draw (19) to (17);
		\draw (17) to (18);
		\draw (20) to (18);
		\draw (18) to (21);
		\draw (21) to (22);
		\draw (21) to (23);
		\draw (23) to (24);
		\draw (23) to (25);
		\draw (25) to (26);
		\draw (25) to (27);
		\draw (34) to (35);
		\draw (37) to (35);
		\draw (35) to (36);
		\draw (38) to (36);
		\draw (36) to (39);
		\draw (39) to (40);
		\draw (39) to (41);
		\draw (41) to (42);
		\draw (41) to (43);
		\draw (43) to (44);
		\draw (43) to (45);
		\draw (51) to (52);
		\draw (54) to (52);
		\draw (52) to (53);
		\draw (55) to (53);
		\draw (53) to (56);
		\draw (56) to (57);
		\draw (56) to (58);
		\draw (58) to (59);
		\draw (58) to (60);
		\draw (60) to (61);
		\draw (60) to (62);
\end{tikzpicture}

    \caption{All four of the above trees $T_1, T_2, T_3, T_4$ are $\pi$-maximising trees for the segment sequence $(25, 11, 10, 4, 1,1,1,1,1,1,1)$. Each has two median vertices in the segment of length 25, and each has proximity $\frac{803}{57} \approx 14.1$. Segments are labelled by length --- the length 1 segments are unlabelled.}
    \label{fig_4_equivalent_trees}
\end{figure}
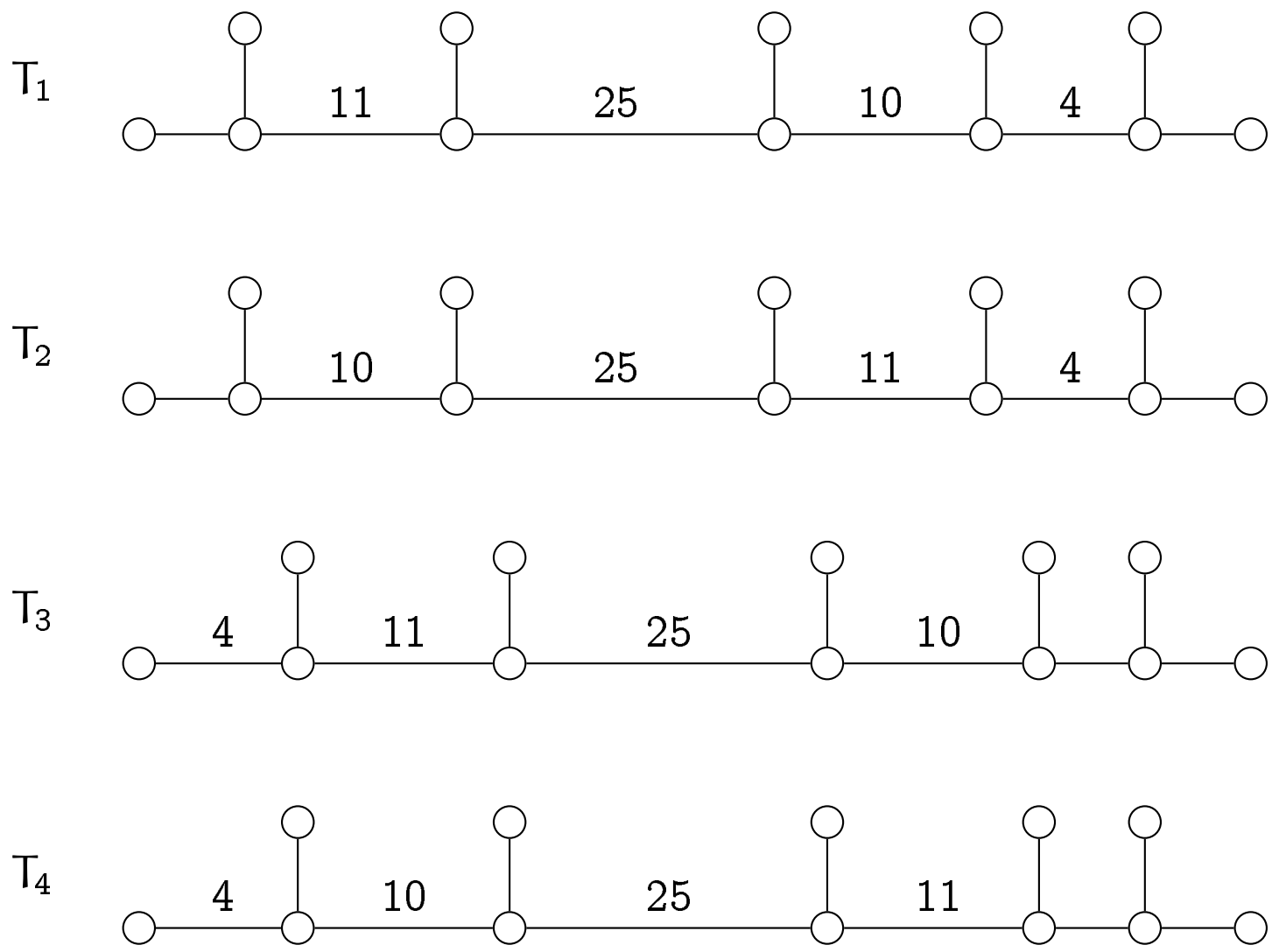

Theorem \ref{thm:pmax_main_theorem} follows immediately from the lemmas that make up the rest of this section. 
Theorem \ref{thm:pmax_main_theorem} provides a detailed, but still incomplete description of any tree that maximizes proximity with given segment sequence $D$. 
The incompleteness arises from the following challenges: the exact arrangement of segments that maximizes proximity is very sensitive to small changes in $D$, and is not even unique! 
For example, Figure \ref{fig_4_equivalent_trees} shows four non-isomorphic trees have segment sequence $25, 11, 10,4,1,1,1,1,1,1,1$ and maximise proximity.

\begin{lem}
    If $T$ is a $D$-tree maximising proximity, then $\Delta(T) \leq 4$.
    \label{lem:pmax_delta_four}
\end{lem}

\begin{proof}
    Assume to the contrary that $T$ is a $\pi$-maximising $D$-tree, and that $T$ has a vertex $v$ of degree at least 5. 
    Let $n$ be the order of $T$.
    Let $w$ be a median vertex of $T$, and $B_1, B_2, \dots, B_m$ be the branches at $w$. 
    There are multiple cases to consider.

    \textit{Case 1:} The vertex $v$ is not $w$.\\
    In this case, $v$ is in a branch $B_i$.
    Assume without loss of generality that $v \in B_1$, and let $P: w, \dots, x$ be a longest path of $T$ that starts at $w$, and ends at a leaf $x$ of $B_1$.
    The vertex $v$ has at least 5 branches, $C_1, C_2, \dots, C_k$ with $k\geq 5$.
    Let $C_k$ be the branch containing $w$, and note that $C_i \leq B_1$ for all $i<k$.
    The path $P$ intersects at most two branches at a vertex other than $v$: $C_k$, and possibly one other branch. 
    Thus, we may assume without loss of generality that $(P\cap C_1) - \{v\} = \emptyset$ and $(P\cap C_2) - \{v\} = \emptyset$.
    Create a new $D$-tree $T'$ by pulling $C_1$ and $C_2$ off of $v$, and attaching them to the leaf $x$. 
    Note that the vertex $w$ is still a median vertex of the new tree $T'$ since the branch-weights of $w$ in $T'$ are the same as the branch-weights of $w$ in $T$. 
    Further, we have that $\sigma(w, T') > \sigma(w, T)$, as all the vertices of $C_1-v$ and $C_2-v$ are further from $w$ in $T'$ than they are from $w$ in $T$.
    Thus $\pi(T') > \pi(T)$, contradicting the maximality of $T$.

    \textit{Case 2:} $v = w$, and at least one branch $B_i$ is trivial.\\
    Without loss of generality, suppose that $B_1=S$ is trivial branch (i.e., consists of a single segment). 
    We create a new tree $D$-tree $T'$. 
    Consider the disjoint union of all the branches $B_1 \sqcup B_2 \sqcup \dots \sqcup B_k$. 
    One leaf of $B_1$ is a copy of $v$, and the other is a vertex $u$. 
    Identify the copies of the root $v$ belonging to $B_1, B_2$ and $B_3$, and identify $u$ with the copies of $v$ belonging to $B_4, B_5, \dots, B_k$.
    Note that this process used to create $T'$ from $T$ is the reverse of the operation from 
    Lemma \ref{lem:rmin_star_operation}.
    Thus, by the same argument as in the proof of Lemma \ref{lem:rmin_star_operation}, $\pi(T) < \pi(T')$.

    \textit{Case 3:} $v = w$, and no branch $B_i$ is trivial.\\
    Suppose without loss of generality that $|B_1| \leq |B_2| \leq \dots \leq |B_k|$. 
    Since no branch is trivial, the branch $B_2$ has at least one branch vertex $u$.
    Create a new $D$-tree $T'$ from $T$ by removing the branch $B_1$ from $v$, and attaching $B_1$ to $u$ (thus we are `merging' branches $B_1$ and $B_2$ into a single larger branch). 
    Since $v=w$ is a median vertex of $T$, we have $|B_i| < \frac{n}{2}$ for all $i$.
    Since there are at least 5 branches, we have $|B_1| + |B_2| < \frac{n}{2}$.
    Thus each branch of $v=w$ in $T'$ has weight less than $\frac{n}{2}$, so $w$ is a median vertex of $T'$.
    We obtained $T'$ from $T$ by moving all the vertices of $B_1$ further away from $w$, so $\sigma(w,T') > \sigma(w,T)$, contradicting the maximality of $T$.
\end{proof}

\begin{lem}
    Let $T$ be a $\pi$-maximising $D$-tree, and suppose $w$ is a median vertex of $T$.
    Let $B$ be any branch of $T$ at $w$.
    Then the following hold:
    \begin{itemize}
        \item There is a path in $B$, starting at $w$, that intersects every segment of $B$,
        \item $B$ has at most 1 vertex of degree 4, and it is further from $w$ than any other branch vertex of $B$,
        \item If $d(w) > 2$, then the longest segment in $B$ is incident with $w$.
    \end{itemize}
    \label{lem:pmax_branch_caterpillar}
\end{lem}

\begin{proof}
    Let $T$ and $w$ be as they are in the statement of the lemma, and suppose that $B$ is a branch of $T$ at $w$.

    \textit{Claim 1:} There is a path in $B$, containing $w$, that intersects every segment of $B$.\\
    Assume to the contrary that the claim is false.
    Let $P = w,\dots, x$ be a maximum length path of $T$ that starts at $w$, and ends at a branch vertex $x\in B$.
    Since the claim is false, the subgraph induced by $E(B) - E(P)$ contains at least one non-trivial tree $C$.
    Note that $C$ is a branch at some vertex $u$ of $P$.
    Since $C$ is a non-trivial branch, it consists of multiple segments, one of which, $C_u$, is incident with $u$.
    Create a new $D$-tree $T'$ from $T$ by removing the subtree of $C - E(C_u)$, and attaching this subtree to the vertex $x$.
    Since we have only moved vertices around in the branch $B$, the vertex $w$ is still a median vertex of $T'$ (its branch-weight is unchanged). 
    As the vertices of $C \setminus C_u$ are further from $w$ in $T'$ than they are from $w$ in $T$, we have $\sigma(w,T') > \sigma(w,T)$.
    This contradicts the maximality of $T$ and proves Claim 1.

    \textit{Claim 2:} $B$ has at most one vertex of degree 4.\\
    Assume to the contrary that there exist two degree 4 vertices, $u$ and $v$.
    Let $Q = w,\dots, y$ be a maximum length path of $T$ that starts at $w$, and ends at a leaf $y\in B$.
    Since $B$ is a subdivided caterpillar (per Claim 1) and $u$ has degree 4, $u$ is incident with at least one pendant segment $S_u$.
    Similarly, $v$ is incident with a pendant segment $S_v$.
    Create a new $D$-tree $T'$ from $T$ by removing $S_u$ and $S_v$ from $u$ and $v$, and attaching these two segments to the vertex $y$. 
    The branch-weight of $w$ is the same in $T'$ as it is in $T$, and the vertices of $S_u-u$ and $S_v-v$ are further from $w$ in $T'$ than they are from $w$ in $T$.
    Thus $\pi(T') > \pi(T)$, contradicting maximality of $T$.
    If $u$ is a branch vertex of $B$ with degree 4, and $v$ is a branch vertex of $B$ with degree 3, it is clear that $d(w,u) > d(w,v)$, as the tree $T'$ formed by removing a segment from $u$ and attaching it to $v$ will have greater proximity.

    \textit{Claim 3:} A segment of maximum length among the segments in $B$ is incident with $w$.\\
    Suppose to the contrary that the segment $R$ of $B$ incident with $w$ does not have maximum length.
    Then there exists another segment $S$ of $B$ with $\ell(S) > \ell(R)$. 
    Create a new tree $T'$ from $T$ by swapping $R$ and $S$.
    Note that the branch-weight of $w$ in $T'$ is the same as the branch-weight of $w$ in $T$, so $w$ is a median vertex of $T'$.
    The vertices that appear after both $R$ and $S$ in a path starting at $w$ will remain the same distance from $w$ in $T'$ as they are from $w$ in $T$.
    All other vertices of $B$ are moved a distance $\ell(S) - \ell(R)$ further from $w$ in $T'$ than they are from $w$ in $T$.
    Thus $\pi(T') > \pi(T)$, contradicting maximality of $T$.
\end{proof}

\begin{lem}
    If a tree is a $\pi$-maximising $D$-tree, then it is a subdivided caterpillar.
    \label{lem:pmax_caterpillar}
\end{lem}

\begin{proof}
    Let $T$ be such a tree, $w$ be a median vertex of $T$, and $n$ the order of $T$.
    By Lemma \ref{lem:pmax_branch_caterpillar}, it suffices to show that at most two of the branches of $T$ at $w$ are non-trivial. 
    So assume to the contrary that $A$, $B$ and $C$ are three non-trivial branches of $T$ at $w$.
    Let $u_a$, $u_b$ and $u_c$ be the branch vertices of $A,B$ and $C$ respectively that are furthest from $w$. 
    Suppose, without loss of generality, that $d(w,u_c) \leq d(w,u_b)$ and $d(w,u_c) \leq d(w, u_a)$.
    By Lemma \ref{lem:pmax_delta_four}, the vertices $u_b$ and $u_c$ have degrees in $\{3,4\}$.

    \textit{Claim:} $w$ is the unique median vertex of $T$ --- hence $bw(w) < \frac{n}{2}$.\\
    Suppose to the contrary that $v\neq w $ is also a median vertex of $T$. 
    Then the branch at $v$ containing $w$ does not have a path intersecting all of its segments --- contradicting Lemma \ref{lem:pmax_branch_caterpillar}.

    \textit{Case 1:} The vertex $u_c$ has degree 4.\\
    The vertex $u_c$ is incident with three pendant segments $C_1$, $C_2$ and $C_3$ with lengths $c_1 \leq c_2 \leq c_3$. 
    Let $h = |C| - (c_1 + c_2 + c_3)$, and note that $h \geq c_3$ per part 3 of Lemma \ref{lem:pmax_branch_caterpillar}.
    Observe that 
    \begin{align*}
        n&\geq |A| + |B| + |C| + 1\\
        &\geq |A| + |B| + h + c_1 + c_2 + c_3 + 1
    \end{align*}
    where the last term comes from counting the vertex $w$ itself.

    \textit{Sub-case 1.1:} $d(w,u_c) < d(w,u_a)$ and $d(w,u_c) < d(w,u_b)$.\\
    Observe that either $|A| + c_1 < \frac{n}{2}$, or $|B| + c_1 < \frac{n}{2}$. 
    For if these were both false, then we would have both $(|A| + c_1) + (|B| + c_1) \geq n$, and $(|A| + c_1) + (|B| + c_1) \leq (|A| + c_1) + (|B| + c_2) < n$, which is a contradiction.
    If $|A| + c_1 < \frac{n}{2}$, then we can form a new tree $T'$ by moving the segment $C_1$ from $u_c$ to the vertex $u_a$.
    Since $|A| + c_1 < \frac{n}{2}$, every branch of $w$ in $T'$ still has branch-weight less than $\frac{n}{2}$, so $w$ is the median of $T'$.
    Further, we obtain $\sigma(w, T') > \sigma(w,T)$, contradicting maximality of $T$.
    The proof when $|B| + c_1 < \frac{n}{2}$ follows similarly.

    \textit{Sub-case 1.2:} $d(w,u_c) = d(w,u_b)$ and $d(u_b) = 4$.\\
    Let $B_1$ be a shortest leaf segment incident with $u_b$, of length $b_1$.
    Per part 3 of Lemma \ref{lem:pmax_branch_caterpillar}, we have 
    \[
        b_1, c_1 \leq d(w, u_a) < |A|.
    \]
    Further, since $c_1$ is the shortest segment incident with $u_c$, and $c_1 \leq d(w, u_c)$ (part 3 of Lemma \ref{lem:pmax_branch_caterpillar}), we have that $c_1 \leq \frac{|C|}{4}$. 
    Similarly, $b_1 \leq \frac{|B|}{4}$.
    Since each branch at $w$ has weight less than $\frac{n}{2}$, we have that $c_1, b_1 < \frac{n}{8}$.
    We thus claim that it is impossible for all three of the following inequalities to hold: $|B| + c_1 \geq \frac{n}{2}$, $|C| + b_1 \geq \frac{n}{2}$ and $|A| + c_1 + b_1 \geq \frac{n}{2}$.
    Suppose to the contrary that they all hold. 
    Then 
    \[
        (|B| + c_1) + (|C| + b_1) + (|A| + c_1 + b_1) \geq \frac{3n}{2}
    \]
    However, we have $|A| + |B| + |C| < n$, and $c_1, b_1 < \frac{n}{8}$, so
    \begin{align*}
        (|B| + c_1) + (|C| + b_1) + (|A| + c_1 + b_1) &= |A| + |B| + |C| + 2(c_1 + b_1)\\
        &< n + 2\left(\frac{n}{8} + \frac{n}{8}\right) = \frac{3n}{2},
    \end{align*}
    Since this is a contradiction, at least one of the inequalities does not hold.
    
    If we have $|C| + b_1 < \frac{n}{2}$, then form a new tree $T'$ by removing the segments $C_1$ and $B_1$ from $T$, and re-attaching both segments to the leaf at the end of $C_2$. 
    Since $|C| + b_1 < \frac{n}{2}$, the branch-weight of $w$ in $T'$ is less than $\frac{n}{2}$, so $w$ is the unique median of $T'$. 
    It is clear that the total distance of $w$ is greater in $T'$ than in $T$, contradicting maximality of $T$.
    
    The case where $|B| + c_1 < \frac{n}{2}$ is very similar --- only we construct $T'$ by moving the segments $C_1$ and $B_1$ onto the end of a segment incident with $u_b$, rather than onto the end of $C_2$.

    Finally, if $|A| + c_1 + b_1 < \frac{n}{2}$, form $T'$ by moving segments $C_1$ and $B_1$ to a leaf at the end of a segment incident with $u_a$.

    \textit{Sub-case 1.3:} $d(w,u_c) = d(w,u_b)$ and $d(u_b) = 3$.\\
    Let $B_1$ and $B_2$ denote the two pendant segments incident with $u_b$, and suppose they have lengths $b_1$ and $b_2$ with $b_1\leq b_2$.
    Note, per part 3 of Lemma \ref{lem:pmax_branch_caterpillar}, that the segment of $B$ incident with $w$ is at least as long as $B_2$.
    We claim that at least one of the following holds:
    \begin{enumerate}
        \item $|A| + b_1 < \frac{n}{2}$ and $|C| + b_2 < \frac{n}{2}$, 
        \item $|A| + b_1 + b_2 < \frac{n}{2}$,
        \item $|C| + b_1 + b_2 < \frac{n}{2}$.
    \end{enumerate}
    Suppose to the contrary that all three are false.
    If $|A| + b_1 \geq \frac{n}{2}$, then by summing this with $|C| + b_1 + b_2 \geq \frac{n}{2}$ yields
    \begin{align*}
        (|A| + b_1) + (|C| + b_1 + b_2) \geq n
    \end{align*}
    However, since a longest segment of $B$ is incident with $w$, we have
    \begin{align*}
        (|A| + b_1) + (|C| + b_1 + b_2) &= |A| + |C| + 2b_1 + b_2\\
            &\leq |A| + |C| + |B|\\
            &< n,
    \end{align*}
    which is impossible.
    If $|C| + b_2 \geq \frac{n}{2}$, we obtain a similar contradiction, so at least one of the three statements is true.

    If $|A| + b_1 < \frac{n}{2}$ and $|C| + b_2 < \frac{n}{2}$, then we create a new tree $T'$ by removing the segments $C_1$, $B_1$ and $B_2$ --- and adding $B_1$ to $u_a$, and both $C_1$ and $B_2$ to the end of the leaf vertex on $C_2$. 
    By our usual branch-weight argument, $w$ is a median vertex of $T'$, and $\sigma(w,T') > \sigma(w,T)$.
    If $|C| + b_1 + b_2 < \frac{n}{2}$, create $T'$ by moving $B_1$ and $B_2$ onto the leaf vertex of segment $C_1$ and obtain a contradiction as $\sigma(w,T') > \sigma(w,T)$.
    The case where $|A| + b_1 + b_2 < \frac{n}{2}$ follows similarly to the $|C| + b_1 + b_2 < \frac{n}{2}$ case.

    \textit{Case 2:} The vertex $u_c$ has degree 3.\\
    Let $C_1$, $C_2$ be the two pendant segments incident with $u_c$, and denote their lengths $c_1 \leq c_2$.
    Let $C_w$ be the segment of $C$ incident with $w$, and note that its length $c_w$ satisfies $c_2 \leq c_w$ per Lemma \ref{lem:pmax_branch_caterpillar}.

    \textit{Sub-case 2.1:} $d(w,u_b) < d(w,u_a)$ or $d(w,u_c) < d(w,u_a)$.\\
    We claim that at least one of the following holds:
    \begin{enumerate}
        \item $|A| + c_1 < \frac{n}{2}$ and $|B| + c_2 < \frac{n}{2}$,
        \item $|A| + c_1 + c_2 < \frac{n}{2}$,
        \item $|B| + c_1 + c_2 < \frac{n}{2}$.
    \end{enumerate}
    Suppose to the contrary that none of the above hold.
    If $|A| + c_1 \geq \frac{n}{2}$ and $|B| + c_1 + c_2 \geq \frac{n}{2}$, then we have that 
    \[
        |A| + |B| + 2c_1 + c_2 \geq n.
    \]
    However, we know that 
    \begin{align*}
        |A| + |B| + 2c_1 + c_2 &\leq |A| + |B| + c_1 + c_2 + c_w\\
            &\leq |A| + |B| + |C| < n.
    \end{align*}
    Thus at least one of the statements is true. 

    If $|A| + c_1 < \frac{n}{2}$ and $|B| + c_2 < \frac{n}{2}$, then we create a new tree $T'$ by removing $C_1$ and $C_2$ from $T$, attaching $C_1$ to $u_a$, and attaching $C_2$ to $u_b$. 
    The argument that $w$ is median vertex of $T'$ and that $\sigma(w, T') > \sigma(w,T)$ is much the same as it is in the previous cases. 
    If \item $|A| + c_1 + c_2 < \frac{n}{2}$, move $C_1$ and $C_2$ onto the end of a pendant segment incident with $u_a$; and if \item $|B| + c_1 + c_2 < \frac{n}{2}$, move $C_1$ and $C_2$ onto a pendant segment incident with $u_b$, and obtain a contradiction as in the previous cases.

    \textit{Sub-case 2.2:} $d(w,u_a) = d(w, u_b) = d(w,u_c)$, and at least one of $u_a$ and $u_b$ have degree 4.\\
    Up to relabelling of the branches of $w$, this case is subsumed by Sub-case 1.3.

    \textit{Sub-Case 2.3:} $d(w,u_a) = d(w, u_b) = d(w,u_c)$, $d(u_a) = d(u_b) = 3$, and the branch $C$ has more than three segments.\\
    Since $d(u_c) = 3$, the path from $w$ to $u_c$ consists of at least two segments: $C_w$, and another segment we call $C_3$ of length $c_3$.

    We claim that $c_2 \leq c_3$. 
    Suppose that the contrary that $c_2 > c_3$. 
    Then the tree $T'$ obtained by swapping $C_3$ and $C_2$ has $\pi(T') > \pi(T)$, which contradicts the maximality of $T$.

    We further claim that if $|A| + c_1 + c_2 \geq \frac{n}{2}$, then $|B| + c_1 + c_2 < \frac{n}{2}$.
    Assume otherwise, and note that
    \[
        |A| + |B| + 2c_1 + 2c_2 \geq n.
    \]
    However, since $c_2 \leq c_3 \leq c_w$:
    \begin{align*}
        |A| + |B| + 2c_1 + 2c_2 &\leq |A| + |B| + c_1 + c_2 + c_3 + c_w\\
            &\leq |A| + |B| + |C| < n.
    \end{align*}
    This contradiction proves our further claim. 

    If $|A| + c_1 + c_2 < \frac{n}{2}$, then we create a new tree $T'$ by removing $C_1$ and $C_2$ from $T$, and attaching them to the end of a pendant segment incident with $u_a$. 
    If $|A| + c_1 + c_2 \geq \frac{n}{2}$, then $|B| + c_1 + c_2 < \frac{n}{2}$, and we create $T'$ by moving $C_1$ and $C_2$ onto the end of a pendant segment incident with $u_b$.
    In both cases, we quickly see that $\pi(T') > \pi(T)$.

    \textit{Sub-case 2.4:} $d(w,u_a) = d(w, u_b) = d(w,u_c)$, $d(u_a) = d(u_b) = 3$, and the branch $C$ has exactly the three segments $C_1$, $C_2$ and $C_w$.\\

    For brevity, let $d = d(w,u_a)$.
    Let $A_1, A_2$ be the pendant segments incident with $u_a$, $B_1, B_2$ the pendant segments incident with $u_b$, and $C_1,C_2$ those incident with $u_c$. 
    Further, let $A_w$, $B_w$ and $C_w$ be the segments of $A$, $B$ and $C$ respectively that are incident with $w$.
    Denote by $a_x$ the length of the segment $A_x$, $x\in \{1,2,w\}$, and denote the lengths $b_x$ and $c_x$ similarly (in particular, $c_w = d$).
    Suppose without loss of generality that $b_w \leq a_w$, and note that $a_w \leq d = c_w$.

    \textit{Claim 1:} $|C| + b_1 + b_2 > \frac{n}{2}$; $|A| + c_1 + c_2 > \frac{n}{2}$; and $|B| + a_1 + a_2 > \frac{n}{2}$.\\
    Assume that $|C| + b_1 + b_2 \leq \frac{n}{2}$. 
    Then we create a new tree $T'$ by removing $B_1, B_2$ from $u_b$ and attaching both to the leaf vertex of $C_1$. 
    Since the branchweight of $w$ in $T'$ is at most $\frac{n}{2}$, $w$ is a median vertex of $T'$. 
    Further, it is clear that $\sigma(w,T') > \sigma(w,T)$, contradicting maximality of $T$.
    The proofs that $|A| + b_1 + b_2 > \frac{n}{2}$ and $|B| + a_1 + a_2 > \frac{n}{2}$ are similar.
    \hfill$\diamond$

    Create a new tree $T'$ from $T$ by removing $C_1$ and $C_2$ from $u_c$, and attaching these two segments to the leaf vertex on the segment $A_1$. 
    Let $w'$ denote the median vertex of $T'$ --- in the event that $T'$ has two median vertices, let $w'$ be the one minimising the distance $d(w,w')$ (since $w'$ is a vertex in $T-(C_1\cup C_2)$, this distance will be the same in $T$ and $T'$).
    Let $i = d(w,w')$. 
    Per Claim 1, $|A| + c_1 + c_2 > \frac{n}{2}$, so the branchweight $bw_{T'}(w) \geq \frac{n}{2}$. 
    Thus, $w\neq w'$, so $i>0$.

    \textit{Claim 2:} $w'$ is a vertex of the segment $A_w$.\\
    Let $P : w = x_0, x_1, \dots, u_a = x_{k-1}$ be the unique path from $w$ to $u_a$ (this path is the same in $T$ and $T'$).  
    Since $w'$ minimises branchweight in $T'$, and $w$ minimises branchweight in $T$, it is clear --- considering the weights of $C_1, C_2$ and $A_1$ --- that $w'$ must lie on $P-w$. 
    Suppose to the contrary that $w'$ lies on $P$, but not $A_w$. 
    Then one branch of $w'$ in $T$ contains $A_w$, $C_w$, and all of $B$. 
    Note that ${a_1, a_2 \leq a_w \leq c_w}$.
    Thus
    \begin{align*}
        bw_{T'}(w') &\geq |B| + c_w + a_w\\ 
            &\geq |B| + a_1 + a_2
    \end{align*}
    However, by Claim 1, $|B| + a_1 + a_2 > \frac{n}{2}$, contradicting that $w'$ is a median vertex of $T'$.
    \hfill$\diamond$

    \textit{Claim 3:} $a_1 \geq i$.\\
    Assume to the contrary that $i>a_1$. 
    Note that $c_w = d$.
    Then, by considering the branch of $T'$ at $w'$ containing $B$, we get that
    \begin{align*}
        bw_{T'}(w') &\geq |B| + c_w + i\\
            &> |B| + a_2 + a_1.
    \end{align*}
    Per Claim 1, $|B| + a_2 + a_1 > \frac{n}{2}$, contradicting that $w'$ is a median of $T'$.
    \hfill$\diamond$

    \textit{Claim 4:} $|B| + c_w - |A| \geq 0$.\\
    Assume the Claim is false, and obtain $|B| < |A| - c_w$. 
    Per Claim 1, we have 
    \[
        \frac{n}{2} < |B| + c_1 + c_2 < |A| - c_w + c_1 + c_2
    \]
    Thus, noting that $c_1, c_2 \leq c_w$, we have
    \begin{align}
        n &< (|B| + c_1 + c_2) + (|A| - c_w + c_1 + c_2)\nonumber\\
            &= |A| + |B| + (2c_1 + 2c_2 - c_w)\nonumber\\
            &\leq |A| + |B| + |C|
        \label{eqn:A_B_cw}
    \end{align}
    Therefore $n< |A| + |B| + |C| = n-1 < n$, which is a contradiction.
    \hfill$\diamond$

    To complete the proof, it suffices to show that $\sigma(w', T') - \sigma(w,T) > 0$. 
    Note that $\sigma(w', T') - \sigma(w,T) = (\sigma(w',T') - \sigma(w,T')) + (\sigma(w,T') - \sigma(w,T))$.
    By construction of $T'$ from $T$, we obtain the following expression for $\sigma(w,T') - \sigma(w,T)$ by considering the contribution of the vertices in $C_1$ and $C_2$ to $\sigma(w,T')$, and subtracting from it their contribution to $\sigma(w,T)$. 
    For any vertex $x$ not in $C_1$ or $C_2$, the contribution of $x$ to $\sigma(w,T') - \sigma(w,T)$ is $0$.
    \begin{align}
        \sigma(w,T') - \sigma(w,T) &= \sum_{j= d+a_1+1}^{d+a_1+c_1}j + \sum_{j= d+a_1+1}^{d+a_1+c_2}j - \sum_{j= d + 1}^{d+c_1}j - \sum_{j= d+1}^{d+c_2}j
            = a_1(c_1 + c_2)
        \label{eqn:first_diff}
    \end{align}
    Observing that the vertices of $A\cup C_1 \cup C_2$ are $i$ closer to $w'$ than $w$ (apart from the $i$ vertices between $w'$ and $w$), and that the vertices of $B$ and $C_w$ are $i$ further from $w'$ than $w$, we have:
    \begin{align}
        \sigma(w',T') - \sigma(w,T') &= i(|B| + c_w) - i(|A| + c_1 + c_2 - i)\nonumber\\
            &= i(|B| + c_w - |A|) +i^2 - i(c_1 + c_2)
        \label{eqn:second_diff}
    \end{align}
    Per Claim 3, we have that $(a_1-i)(c_1 + c_2) \geq 0$, so summing equations \ref{eqn:first_diff} and \ref{eqn:second_diff} yields:
    \begin{align*}
        \sigma(w', T') - \sigma(w,T) &\geq i(|B| + c_w - |A|) + i^2
    \end{align*}
    Per Claim 4, $|B| + c_w - |A| \geq 0$, and we know $i>0$.
    Therefore $\sigma(w', T') - \sigma(w,T) > 0$, contradicting that $T$ maximises proximity and completing the proof.
\end{proof}

\begin{lem}
\label{lem:pmax_segment_orders}
    Suppose $T$ is a $\pi$-maximising $D$-tree with median vertex $w$, and that $B$ is a branch of $T$ at $w$.
    Suppose that $S_i, S_j$ are two segments of $T$, fully contained in $B$, with lengths $\ell_i$ and $\ell_j$ respectively. 
    Among paths of $B$ that start at $w$ and meets every segment of $B$, let $P$ be one of maximum length.
    All the following hold:
    \begin{enumerate}
        \item If $S_i$ and $S_j$ are both leaf segments, then $d(w, S_i) < d(w, S_j)$ implies $\ell_i \leq \ell_j$.
        \item If $S_i$ and $S_j$ are segments of $P$, then $d(w, S_i) < d(w, S_j)$ implies $\ell_i \geq \ell_j$.
        \item If $S_i$ is a segment of $P$ and $S_j$ is a leaf segment, then $\ell_i \geq \ell_j$.
    \end{enumerate}
\end{lem}

\begin{proof}
    To prove (1), suppose to the contrary that $S_i, S_j$ are leaf segments with where $S_i$ is closer to $w$ than $S_j$, but $\ell_i > \ell_j$.
    Create $T'$ by removing $(\ell_i - \ell_j)$ vertices from $S_i$ and appending these vertices to the end of $S_j$. I.e., up to renaming the vertices, we have swapped the segments $S_i$ and $S_j$. 
    The branch-weight at $w$ is the same in $T'$ as in $T$, and the $(\ell_i - \ell_j)$ vertices that were moved are further from $w$ in $T'$ than they are in $T$, so $\pi(T') > \pi(T)$, contradicting our choice of $T$.

    The argument to prove (2) is similar: assume to the contrary that $S_i$ and $S_j$ are segments of $P$ such that $d(w, S_i) < d(w, S_j)$ and $\ell_i < \ell_j$.
    We again remove $\ell_j - \ell_i$ vertices from $S_j$ and insert them into $S_i$ to create a new tree $T'$ (up to the vertex labels, this is just swapping $S_i$ and $S_j$). 
    Note that the length of the path $P$ is unaffected by this change, and so the vertices of $P$ contribute the same amount to $\sigma(w,T)$ and $\sigma(w,T')$. 
    However, the leaf segments incident with $P$ between $S_j$ and $S_i$ (of which there is at least one) are $\ell_j - \ell_i$ further from $w$ in $T'$ than they are from $w$ in $T$. 
    Thus $\pi(T) < \pi(T')$, a contradiction.

    We now prove (3). 
    Let $v$ denote the branch vertex of $B$ furthest from $w$ (I.e., $v$ is the last branch vertex of the path $P$). 
    By (1), a longest leaf segment $S_l$ of $B$ is incident with $v$, and is contained in $P$. 
    By (2), any other segment $S_p$ of $P$ is at least as long as the last segment $S_l$ of $P$, completing the proof.
\end{proof}

\begin{lem}
\label{lem:pmax_longest_segment}
    Suppose $T$ is a $\pi$-maximising $D$-tree and that $w$ is a median vertex of $T$.
    Then $w$ lies in a longest segment of $T$.
\end{lem}

\begin{proof}
    Let $S_1, \dots, S_k$ be the segments of $T$ with lengths $\ell_1 \geq \dots \geq \ell_k$ and suppose to the contrary that $w$ is a median vertex that is in a segment $S_j$ such that $\ell_j < \ell_1$.
    Denote the vertices of $S_j$ by $S_j = u_0, \dots, u_{\ell_j}$.

    Per Lemma \ref{lem:pmax_branch_caterpillar}, $T$ is a caterpillar.
    If $w$ is an end-vertex of the segment $S_j$, then by Lemma \ref{lem:pmax_segment_orders}, $w$ is also the end-vertex of a longest segment of each branch at $w$, and is thus in a longest segment of $T$. 
    Thus we assume that $w = u_i$ where $1\leq i \leq \ell_j - 1$.

    The tree $T$ has two branches at $w$, which we call $A$ and $B$. 
    Suppose without loss of generality that $S_1$ is in branch $A$. 
    Per Lemma \ref{lem:pmax_segment_orders}, we may assume $S_1$ is incident with $u_0$.
    Denote $S_1 = v_0, \dots, v_{\ell_1}$, where $v_{\ell_1} = u_0$, and note that there is a leaf segment (or two) incident with $u_0$.
    Create a new tree $T'$ by removing the leaf segments from $u_0$, and attaching them to the vertex $v_{\ell_j}$ (I.e., up to renaming the vertices, we have swapped $S_1$ and $S_j$).
    As we have only moved a leaf segment (or two segments) further from $w$ inside the branch $A$, the weights of the branches $A$ and $B$ are unchanged, so $w$ is a vertex of minimum proximity in $T'$. 
    Further, since the leaf segments are a distance $\ell_1 - \ell_j > 0$ further from $w$ in $T'$ than they from $w$ in $T$, and no other vertices were moved, we conclude that $\pi(T') > \pi(T)$. 
\end{proof}

\medskip

\section{Conclusion}

In this paper, we have studied the proximity and remoteness of trees with
a prescribed segment sequence. We have characterised the unique tree that
minimises proximity and the trees that maximise remoteness. For minimum
remoteness, we have established sufficient conditions under which a
minimising tree has two branching vertices. Our results suggest that the
structure of a remoteness-minimising tree depends on the relative lengths
of the segments. For maximum proximity, we have obtained a detailed set
of necessary conditions. However, these conditions do not give a complete
characterisation, since different non-isomorphic trees can maximise
proximity for the same segment sequence.

An open question concerns the structure of trees that minimise remoteness.
We conjecture that, for every segment sequence
$D=(\ell_1,\ell_2,\ldots,\ell_k)$, a $D$-tree with minimum remoteness is
either starlike or a $(\Delta,2)$-$D$-tree. In particular, it remains to
determine whether a tree with three or more branching vertices can ever
minimise remoteness. Another direction for future research is to seek a
complete characterisation of trees that maximise proximity. The sensitivity
of the structure of proximity-maximising trees to changes in the segment
lengths, together with the possibility of multiple non-isomorphic extremal
trees, makes this a challenging problem. Further investigation of these
questions may lead to a more complete understanding of the extremal
proximity and remoteness of trees with prescribed segment sequences.

\end{document}